\documentclass[10pt,a4paper]{article}

\usepackage{graphics}
\usepackage{graphicx,epsfig}
\usepackage{epsfig}
\usepackage{graphicx}
\usepackage{a4wide}
\usepackage{amsmath}
\usepackage{amsfonts}
\usepackage{cite}
\usepackage{enumerate}
\usepackage{comment}
\usepackage{lipsum}

\begin{document}

\newcommand{\lap}{\Delta}
\newcommand{\be}{\begin{equation}}
\newcommand{\ee}{\end{equation}}
\newcommand{\bee}{\begin{equation*}}
\newcommand{\eee}{\end{equation*}}
\newcommand{\bea}{\begin{eqnarray}}
\newcommand{\eea}{\end{eqnarray}}
\newcommand{\bess}{\begin{eqnarray*}}
\newcommand{\eess}{\end{eqnarray*}}
\numberwithin{equation}{section}
\topmargin=-0.4truecm
\textwidth=15truecm
\oddsidemargin=0.5truecm
\evensidemargin=0.5truecm
\textheight=24.3truecm
\footskip=1truecm



\newcommand{\qed}{\hfill \rule{2mm}{2mm}}
\newcommand{\pf}{{\bf Proof:}}
\newtheorem{definition}{Definition}
\newtheorem{Thm}{Theorem}[section]
\newtheorem{Lem}{Lemma}[section]
\newtheorem{proposition}{Proposition}[section]
\newtheorem{remark}{Remark}[section]
\newtheorem{corollary}{Corollary}[section]
\newtheorem{example}{Example}[section]
\title{\bf 
Integral representations of solutions to \\ Poisson equations
}
\author{$^{1}$Aye Chan May and $^{2^{*}}$Adisak Seesanea\\
$^{1,2^{*}}$Sirindhorn International Institute of Technology \\
Thammasat University, Thailand\\
e-mail : $^{1}$m6422040748@g.siit.tu.ac.th / $^{2^{*}}$adisak.see@siit.tu.ac.th\\
$^{*}$Corresponding author}
\date{}

\maketitle

\thispagestyle{empty}

\begin{abstract}
\noindent We give a constructive approach for the study of integral representations of classical solutions to Poisson equations under some integrability conditions on data functions. 
\end{abstract}

\section{Introduction}
The fundamental solution of the Laplace equation, traditionally defined for $\mathbb{R}^n$, has been a focal point in mathematical research due to its fundamental role in partial differential equations. Historically, the exploration of this solution has primarily centered around functions within the class $\mathcal{C}_0^2(\mathbb{R}^n)$. The roots of this research lie in the classical treatment of Laplace's equation and the associated Poisson problem, where the fundamental solution is a foundation.

We study the solutions of the Poisson equation involving the Laplace operator 
\be\label{poisson}
\Delta u=f\quad\text{in}\quad\mathbb{R}^n
\ee
under the large classes of the functions $f$. Thus a solution $u$ of the equation \eqref{poisson} having compact support can always be recovered from the right-hand side $f$ in terms of a convolution of $f$ with a fundamental solution $E_n(x)$.

The  {\it fundamental solution for the Laplacian} in $\mathbb{R}^n$ which is defined by
\[
E_{n}(x) := 
\begin{cases}
   \frac{-1}{(n-2)\omega_{n-1}}\frac{1}{\left\vert x \right\vert^{n-2}} \qquad \text{for}\;\; x \in \mathbb{R}^n \setminus \lbrace 0 \rbrace \;\; \text{if} \;\; n \geq 3, \\
  \frac{1}{2\pi}\ln |x| \,\;\quad\qquad\qquad \text{for}\;\; x \in \mathbb{R}^n \setminus \lbrace 0 \rbrace \;\; \text{if} \;\; n =2,
\end{cases}
\]
where $\omega_{n-1}$ is the surface area of the unit sphere in $\mathbb{R}^n$. 
It is well-known that if $f \in \mathcal{C}^2_{0}(\mathbb{R}^n)$, then the function 
\be\label{solu}
u(x):=  (E_n \ast f)(x) =\int_{\mathbb{R}^n} E_{n}(x-y)f(y) \;dy \qquad\qquad \text{for}\;\;x \in \mathbb{R}^n.
\ee
is well-defined in $\mathbb{R}^n$, and belongs to $\mathcal{C}^2(\mathbb{R}^n)$, see \cite{Ev}. 

The seminal works of De Giorgi \cite{De} and Nash \cite{Nash} are widely recognized for their significant contributions to establishing the H\"{o}lder regularity of solutions under specific assumptions on the function $f$. It has been observed that the function $u$ belongs to the class $C^{2}(\Omega)$ under the condition that the function $f$ is bounded and locally H\"{o}lder continuous in a bounded domain $\Omega\subset\mathbb{R}^{n}$, as stated in \cite[Lemma 4.2]{GT}. On the other hand, the authors in \cite{AM} proved that the solution $u$ of \eqref{poisson} is in the class of $C^{1}(\Omega)$ and H\"{o}lder continuous if $f$ is bounded and continuous.
Further related results can be found in \cite{EO},\cite{MA},\cite{LW}, \cite{DM}, \cite{STT},\cite{MMS}, and the literature cited therein.

Our work is mainly motivated by the contributions of De Giorgi \cite{De}, Nash \cite{Nash}, Gilbarg and Trudinger \cite{GT}, and Evan \cite{Ev}. In this paper, we introduce the convenient conditions of $f$ which is certainly integrable and continuously differentiable. This allows us to establish regularity results to cover functions $u$ defined over the entire space $\mathbb{R}^{n}$. As a consequence, we establish the existence and uniqueness of bounded solutions to Dirichlet problem for Poisson equation with the source $f$ belongs to Lorentz spaces.




The following theorem illustrates the well-defined nature of solutions throughout various function spaces, providing perspectives on existence, continuity, and higher-order differentiability practically everywhere.
\begin{Thm}\label{mainthm}
Let $f: \mathbb{R}^n \rightarrow \mathbb{R}$ be a Lebesgue measurable function. Define 
\[
u(x):= \int_{\mathbb{R}^n} E_{n}(x-y)f(y) \;dy \qquad\qquad \text{for}\;\;x \in \mathbb{R}^n.
\]
 If $\frac{f(y)}{1+|y|^{n-2}} \in L^1(\mathbb{R}^n)$, $\frac{\nabla f(y)}{1+|y|^{n-1}} \in L^1(\mathbb{R}^n)$ and $f \in \mathcal{C}^{1}(\mathbb{R}^n)$, then $u$ is well-defined in $\mathbb{R}^n$ and belongs to $\mathcal{C}^{2}(\mathbb{R}^n)$. Moreover, for each $i,j \in \lbrace 1,...,n\rbrace $ and $x \in \mathbb{R}^n$ we have 
 \[
 \partial_{i}\partial_{j} u (x) = \int_{\mathbb{R}^n} \partial_{j} E_{n}(x-y) \partial_{i}f (y)\;dy.
 \]
 Furthermore, $\Delta u = f $ in $\mathbb{R}^n$.
 \end{Thm}
 \noindent{\bf Organization of paper}\\

In Section \ref{sec2}, we give some notations and well-known results related to our work. The proof of the main theorem and necessary lemmas are discussed in Section \ref{1-sec-2}. Moreover, we also give an application of the main theorem to the Poisson problem.

 \section{Preliminaries} \label{sec2}
 Throughout this paper, let $\mathbb{R}^n$ be Euclidean space with $n\geq 3.$ We denote $\mathcal{C}_0^2(\mathbb{R}^n)$ is the space of smooth functions with compact support in $\mathbb{R}^n$. Given a Lebesgue measurable set $E\subset\mathbb{R}^n$, the integral of a measurable function $f: E\to \mathbb{R}$ with respect to Lebesgue measure can be defined as
 \[
 \int_E fdx = \int_{\mathbb{R}^n}  \tilde{f}dx,\quad\text{where}\quad\tilde{f}(x) = \begin{cases} f(x)\quad\text{if}\;\; x\in E,\\ 0\qquad \;\; \text{if}\;\;x\not\in E.
 \end{cases}
 \]
 The space of all Lebesgue measurable functions $f:E \to \mathbb{R}$ such that
$
 \|f\|_{L^p(E)}:=\big(\int_E |f|^p dx\big)^{\frac{1}{p}}<\infty,
$
 for $1\leq p<\infty$. Given an open set $\Omega\subset\mathbb{R}^n$, we denote $L^p_{loc}(\Omega)$ as the space of all measurable functions $f:\Omega\to \mathbb{R}^n$ such that $f\in L^p(U)$ for every bounded open set $U$ whose closure is contained in $\Omega$.
 
 \begin{Thm}[\textbf{Mean-value formulas for Laplace eqaution}]\label{mean}
 If $u\in\mathcal{C}^2(\mathbb{R}^n)$ is harmonic, then
 \[
 u(x)=\frac{1}{\omega_{n-1} r^{n-1}}\int_{\partial B(x,r)} u\;dS=\frac{1}{\omega_n r^n}\int_{B(x,r)} u\;dy,
 \]
for each ball $B(x,r)\subset\mathbb{R}^n$, $\omega_{n-1}$ is surface area of the unit sphere and $\omega_n$ is volume of unit ball in $\mathbb{R}^n$.
\end{Thm}
A complete proof of the following regularity result can be found in \cite{Ev}.
\begin{Thm}
 If $u\in \mathcal{C}(\mathbb{R}^n)$ satisfies the mean-value property for each ball $B(x,r)\subset \mathbb{R}^n$, then
$
 u\in \mathcal{C}^\infty(\mathbb{R}^n).
$
 \end{Thm}

\section{Fundamental Solutions}\label{1-sec-2}
In this section, we begin by giving the important lemmas for our main theorem.
\begin{Lem}\label{lemma1}
Let $f: \mathbb{R}^n \rightarrow \mathbb{R}$ be a Lebesgue measurable function. Define 
\[
u(x):= \int_{\mathbb{R}^n} E_{n}(x-y)f(y) \;dy \qquad\qquad \text{for}\;\;x \in \mathbb{R}^n.
\] If  $\frac{f(y)}{1+|y|^{n-2}} \in L^1(\mathbb{R}^n) $, then $u$ is well-defined almost everywhere in $\mathbb{R}^n$.
\end{Lem}
\noindent{\bf Proof.}
Let $R>0$ and $x \in B(0,\frac{R}{2})$. Write
\[
\begin{split}
\int_{\mathbb{R}^n} |E_{n}(x-y)||f(y)| \;dy 
&= \frac{1}{(n-2)\omega_{n-1}}\int_{\mathbb{R}^n} \frac{|f(y)|}{|x-y|^{n-2}} \;dy \\
&= \frac{1}{(n-2)\omega_{n-1}} \left[ \int_{B(0,R)} \frac{|f(y)|}{|x-y|^{n-2}} \;dy + \int_{\mathbb{R}^n \setminus B(0,R)} \frac{|f(y)|}{|x-y|^{n-2}} \;dy \right] \\
&= \frac{1}{(n-2)\omega_{n-1}} \left[\; I_{R}(x) + II_{R}(x)\; \right],
\end{split}
\]
where 
\[
I_{R}(x):=\int_{B(0,R)} \frac{|f(y)|}{|x-y|^{n-2}} \;dy \qquad\text{and} \qquad II_{R}(x):=\int_{\mathbb{R}^n \setminus B(0,R)} \frac{|f(y)|}{|x-y|^{n-2}} \;dy.
\]
Now, we estimate the second term $II_{R}(x)$. If $y \in \mathbb{R}^n \setminus B(0,R)$, then 
\be\label{ineq1}
|x-y| \geq |y|-|x| > R - \frac{R}{2} = \frac{R}{2},\qquad\text{and}
\ee
\be\label{ineq2}
\qquad  |x-y| \geq |y|-|x| = \frac{|y|}{2} + \frac{|y|}{2} - |x|  > \frac{|y|}{2}
\qquad \qquad \qquad \qquad \;\;
\ee
which implies that 
\[
|x-y|^{n-2} > \frac{1}{2} \left[ \left(\frac{R}{2}\right)^{n-2} + \left( \frac{|y|}{2} \right)^{n-2} \right] 
 \geq  c_{1} \left(1+|y|^{n-2} \right),
\]
where $c_{1}=c_{1}(n,R):=\frac{1}{2^{n-1}}\;\text{min} \lbrace R^{n-2},1 \rbrace$. Thus, 
\[
II_{R}(x) = \int_{\mathbb{R}^n \setminus B(0,R)} \frac{|f(y)|}{|x-y|^{n-2}} \;dy \leq  \frac{1}{c_{1}}\int_{\mathbb{R}^n \setminus B(0,R)} \frac{|f(y)|}{1+ |y|^{n-2}} \;dy \leq \frac{1}{c_{1}}\int_{\mathbb{R}^n} \frac{|f(y)|}{1+ |y|^{n-2}} \;dy < +\infty.
\]
This shows that $II_{R}(x) < + \infty $ for all $x \in B(0,\frac{R}{2})$. Next, we estimate the first term $I_{R}(x)$. Note that
\[
\begin{split}
\int_{B(0,\frac{R}{2})} I_{R}(x) \;dx 
&= \int_{B(0,\frac{R}{2})} \int_{B(0,R)} \frac{|f(y)|}{|x-y|^{n-2}} \;dy\;dx \\
&= \int_{B(0,R)} |f(y)| \int_{B(0,\frac{R}{2})}  \frac{1}{|x-y|^{n-2}} \;dx\;dy \\
&\leq \int_{B(0,R)} |f(y)| \int_{B(y,\frac{3R}{2})}  \frac{1}{|x-y|^{n-2}} \;dx\;dy \\
&= \int_{B(0,R)} \frac{|f(y)|}{1+|y|^{n-2}} \left( 1+|y|^{n-2}\right) \int_{B(y,\frac{3R}{2})}  \frac{1}{|x-y|^{n-2}} \;dx\;dy \\
&\leq \left( 1+ R^{n-2}\right)\int_{B(0,R)} \frac{|f(y)|}{1+|y|^{n-2}}  \int_{B(0,\frac{3R}{2})}  \frac{1}{|z|^{n-2}} \;dz\;dy \\
&= \left( 1+ R^{n-2}\right) \left( \int_{B(0,R)} \frac{|f(y)|}{1+|y|^{n-2}}  \;dy \right) \left( \int_{B(0,\frac{3R}{2})}   \frac{1}{|z|^{n-2}} \;dz \right) \\
&\leq \left( 1+ R^{n-2}\right) \left( \int_{\mathbb{R}^n} \frac{|f(y)|}{1+|y|^{n-2}}  \;dy \right) \left( \int_{0}^{\frac{3R}{2}}\int_{S^{n-1}}   \frac{\rho^{n-1}}{\rho^{n-2}} \;d\sigma(\omega)\;d\rho \right) \\
&= \left( 1+ R^{n-2}\right)\omega_{n-1} \left( \int_{\mathbb{R}^n} \frac{|f(y)|}{1+|y|^{n-2}}  \;dy \right) \left( \int_{0}^{\frac{3R}{2}}   \rho\;d\rho \right) \\
&= \frac{9R^{2}\left( 1+ R^{n-2}\right)\omega_{n-1}}{8}\left( \int_{\mathbb{R}^n} \frac{|f(y)|}{1+|y|^{n-2}}  \;dy \right) < +\infty.
\end{split}
\]
This implies that $I_{R}(x) < +\infty$ for almost every $x \in B(0,\frac{R}{2})$. Hence,
\[
\int_{\mathbb{R}^n} |E_{n}(x-y)||f(y)| \;dy < +\infty \qquad\qquad \text{for almost every}\;\;x \in B(0,\frac{R}{2}).
\] 
Since $R>0$ was arbitrary, we have shown, in particular, that for every $j \in \mathbb{N}$ there is a null set $N_{j} \subset B(0,\frac{j}{2})$ such that 
\[
\int_{\mathbb{R}^n} |E_{n}(x-y)||f(y)| \;dy < +\infty \qquad\qquad \text{for all}\;x \in B(0,\frac{j}{2})\setminus N_{j}.
\]
Setting $N:=\cup_{j=1}^{\infty}N_{j}$, which is a null set in $\mathbb{R}^n$ so that 
\[
\int_{\mathbb{R}^n} |E_{n}(x-y)||f(y)| \;dy < +\infty \qquad\qquad \text{for all}\;x \in \mathbb{R}^n \setminus N.
\]
This proves that $u$ is well-defined almost everywhere in $\mathbb{R}^n$.

\begin{Lem}\label{lemma2}
Let $f: \mathbb{R}^n \rightarrow \mathbb{R}$ be a Lebesgue measurable function. Define 
\[
u(x):= \int_{\mathbb{R}^n} E_{n}(x-y)f(y) \;dy \qquad\qquad \text{for}\;\;x \in \mathbb{R}^n.
\]Consider  $\frac{f(y)}{1+|y|^{n-2}} \in L^1(\mathbb{R}^n) $ and $f \in L^{r}_{loc}(\mathbb{R}^n)$ for some $r \in (\frac{n}{2},\infty] $. Then $u$ is well-defined and continuous in $\mathbb{R}^n$.
\end{Lem}
\noindent{\bf Proof.} In view of Lemma \ref{lemma1}, we see that once we have established that $I_{R}(\cdot) < +\infty$ in each ball $B(0,\frac{R}{2})$, then $u$ is well-defined everywhere on $\mathbb{R}^n$. Let $R>0$ and $x \in B(0,\frac{R}{2})$. Applying H{\"o}lder inequality, we get 
\[
I_{R}(x) \leq \left\lVert f \right\rVert_{L^{r}(B(0,R))} \left\lVert \frac{1}{|x-y|^{n-2}}\right\rVert_{L^{r'}(B(0,R))} \quad\text{where}\quad \frac{1}{r} + \frac{1}{r'} =1.
\]
Since $f \in L^{r}_{loc}(\mathbb{R}^n)$, it is left to show that $\left\lVert \frac{1}{|x-y|^{n-2}}\right\rVert_{L^{r'}(B(0,R))}$  is finite. To see this,
\[
\left\lVert \frac{1}{|x-y|^{n-2}}\right\rVert_{L^{r'}(B(0,R))}^{r'} 
= \int_{B(0,R)} \frac{1}{|x-y|^{r'(n-2)}}\;dy 
\leq \int_{B(0,\frac{3R}{2})} \frac{1}{|z|^{r'(n-2)}}\;dz < +\infty,
\]
where the last integral on the right-hand side is finite because $r > \frac{n}{2}$, equivalently, $r'(n-2) < n$. Hence, by the discussion above, $u$ is well-defined in $\mathbb{R}^n$. Next, we show that $u$ is continuous in $\mathbb{R}^n$. Let $x_{0} \in \mathbb{R}^n$, and let a sequence $\lbrace x_k \rbrace_{k \in \mathbb{N}} \subset \mathbb{R}^n$ such that $x_k \rightarrow x_0$ as $k \rightarrow \infty$. We need to show that 
\[
\int_{\mathbb{R}^n} E_{n}(x_{k}-y)f(y) \;dy \longrightarrow 
\int_{\mathbb{R}^n} E_{n}(x_{0}-y)f(y) \;dy \qquad \text{as} \;\; 
k\rightarrow \infty.
\]
For $R>0$ and $k \in \mathbb{N}$, set
\[
I_{k,R}:= \int_{B(0,R)} E_{n}(x_{k}-y)f(y) \;dy, \qquad\;\;
II_{k,R}:=\int_{\mathbb{R}^n \setminus B(0,R)} E_{n}(x_{k}-y)f(y) \;dy
\]
and 
\[
I_{0,R}:= \int_{B(0,R)} E_{n}(x_{0}-y)f(y) \;dy, \qquad\;\;
II_{0,R}:=\int_{\mathbb{R}^n \setminus B(0,R)} E_{n}(x_{0}-y)f(y) \;dy.
\]
Note that 
\[
\int_{\mathbb{R}^n} E_{n}(x_{k}-y)f(y) \;dy = I_{k,R} + II_{k,R}
\qquad
\text{and} 
\qquad
\int_{\mathbb{R}^n} E_{n}(x_{0}-y)f(y) \;dy = I_{0,R} + II_{0,R},
\]
and also note that since $\lbrace x_{k}\rbrace_{k \in \mathbb{N}}$ is convergent, then it is bounded. Pick $M \in (0,\infty)$ such that $\lvert x_k \rvert<M$ for all $k \in \mathbb{N}$.
Let $\epsilon > 0$. First, we claim that 
\[
\text{sup}_{k \in \mathbb{N}} \left| II_{k,R} \right| < \frac{\epsilon}{3} \qquad \text{and} \qquad |II_{0,R}| < \frac{\epsilon}{3} \qquad \text{for}\;\;R\;\; \text{large enough}.
\]
For every $R > 2M,$ $y \in \mathbb{R}^n \setminus B(0,R)$ and $k \in \mathbb{N}$ we have 
\be\label{ineq3}
|x_{k}-y| \geq |y|-|x_{k}| > R-M > M
\ee
and
\be\label{ineq4}
 |x_{k}-y| \geq |y|-|x_{k}| = \frac{|y|}{2} + \frac{|y|}{2} - |x_{k}| > \frac{|y|}{2} + \frac{R}{2} - M > \frac{|y|}{2},
\ee
which implies that 
\[
|x_{k}-y|^{n-2} \geq \frac{1}{2} \left[ M^{n-2} + \left( \frac{|y|}{2} \right)^{n-2} \right] 
 \geq  c_{2} \left(1+|y|^{n-2} \right),
\]
where $c_{2} = c_{2}(n,M) := \frac{1}{2^{n-1}} \text{min} \lbrace (2M)^{n-2}, 1 \rbrace $. Thus, for $R>2M$ we have
\[
\text{sup}_{k \in \mathbb{N}} |II_{k,R}| 
\leq \frac{1}{(n-2)\omega_{n-1}} \text{sup}_{k \in \mathbb{N}}  \int_{\mathbb{R}^n \setminus B(0,R)} \frac{|f(y)|}{|x_{k}-y|^{n-2}} \;dy 
\leq  \frac{1}{c_{2}(n-2)\omega_{n-1}}\int_{\mathbb{R}^n \setminus B(0,R)} \frac{|f(y)|}{1+ |y|^{n-2}} \;dy.
\]
Since $ \frac{f(y)}{1+ |y|^{n-2}} \in L^{1}(\mathbb{R}^n) $, then by the Lebesgue Dominated Convergence Theorem, we have
\[
\int_{\mathbb{R}^n \setminus B(0,R)} \frac{|f(y)|}{1+ |y|^{n-2}} \;dy = \int_{\mathbb{R}^n} \frac{|f(y)|}{1+ |y|^{n-2}}\chi_{\mathbb{R}^n \setminus B(0,R)}(y) \;dy  \longrightarrow 0 \qquad \text{as}\;\;R \rightarrow \infty.
\]
So there exists $R_{1} > 0$ such that 
$
\text{sup}_{k \in \mathbb{N}} |II_{k,R}| < \frac{\epsilon}{3} \qquad \text{for}\;\;R \geq R_1.
$
Next, to get our second claim, we apply the same argument as above. Then, for every $R > 2|x_{0}|$ and $y \in \mathbb{R}^n \setminus B(0,R)$,
there exists $R_{2} > 0$ such that 
$
|II_{0,R}| < \frac{\epsilon}{3} \qquad \text{for}\;\;R \geq R_2.
$
Setting $\tilde{R}:= \text{max}\lbrace R_{1}, R_{2}\rbrace$, then 
\[
\text{sup}_{k \in \mathbb{N}} |II_{k,\tilde{R}}| < \frac{\epsilon}{3} \qquad \text{and} \qquad |II_{0,\tilde{R}}| < \frac{\epsilon}{3}.
\]
Finally, we show that 
$
I_{k,\tilde{R}} \longrightarrow I_{0,\tilde{R}} \qquad \text{as}\;\;k \rightarrow \infty,
$
using the Vitali's Convergence Theorem on $B(0,\tilde{R})$. 
For $y \in B(0,\tilde{R})$ and $ k \in \mathbb{N}$, set  
\[
F_{k}(y):= E_{n}(x_{k}-y)f(y) \qquad \text{and} \qquad F_{0}(y):= E_{n}(x_{0}-y)f(y),
\]
which are the kernels of $I_{k,\tilde{R}}$ and $I_{0,\tilde{R}}$, respectively. Clearly that $F_{k} \rightarrow F_{0}$ pointwise almost everywhere on $B(0,\tilde{R})$, as $k \rightarrow \infty$. Now, let us check the uniform integrability of the sequence $\lbrace F_{k}\rbrace_{k \in \mathbb{N}}$. Let $\eta >0$, and let $A \subset B(0,\tilde{R})$. For each $k \in \mathbb{N}$, applying H{\"o}lder's inequality, we have
\[
\begin{split}
\left\vert \int_{A} F_{k}(y) \;dy \right\vert&\leq \frac{1}{(n-2)\omega_{n-1}}\int_{A} \frac{|f(y)|}{|x_{k}-y|^{n-2}} \;dy \leq \frac{1}{(n-2)\omega_{n-1}} \|f\|_{L^{r}(A)} \left\lVert\frac{1}{|x_{k}-y|^{n-2}} \right\rVert_{L^{r'}(A)} \\ 
&\leq \frac{1}{(n-2)\omega_{n-1}} \|f\|_{L^{r}(B(0,\tilde{R}))} \left\lVert \frac{1}{|x_{k}-y|^{n-2}} \right\rVert_{L^{r'}(A)} \leq c_{4} \left\lVert \frac{1}{|x_{k}-y|^{n-2}} \right\rVert_{L^{r'}(A)},
\end{split}
\]
where $c_{4} = c_{4}(n,r, \tilde{R}) := \frac{1}{(n-2)\omega_{n-1}} \|f\|_{L^{r}(B(0,\tilde{R}))}$ which is finite as $f \in L^{r}_{loc}(\mathbb{R}^n)$. Note that since $r > \frac{n}{2}$, equivalently, $r'(n-2) < n$, that is, $\frac{n}{r'(n-2)} > 1$. Let $p \in \left( 1,\frac{n}{r'(n-2)} \right)$ and $k \in \mathbb{N}$, applying H{\"o}lder's inequality again, we have
\[
\begin{split}
\left\lVert \frac{1}{|x_{k}-y|^{n-2}} \right\rVert_{L^{r'}(A)}^{r'} 
&= \int_{A} \frac{1}{|x_{k}-y|^{r'(n-2)}}\;dy \\
&= \int_{x_{k}-A} \frac{1}{|z|^{r'(n-2)}}\;dz \\
&= \int_{B(0,\tilde{R}+M)} \frac{1}{|z|^{r'(n-2)}} \chi_{x_{k}-A}(z)\;dz \\
&\leq \left( \int_{B(0,\tilde{R}+M)} \frac{1}{|z|^{pr'(n-2)}}\;dz \right)^{\frac{1}{p}}  \left( \int_{B(0,\tilde{R}+M)} \left( \chi_{x_{k}-A}(z)\right)^{p'} \;dz \right)^{\frac{1}{p'}} = c_{5} \left( \mu \left( A \right) \right)^{\frac{1}{p'}} ,
\end{split}
\]
where $\frac{1}{p} + \frac{1}{p'} = 1$ and $c_{5}=c_{5}(n,\tilde{R},M,r,p):= \left( \int_{B(0,\tilde{R}+M)} \frac{1}{|z|^{pr'(n-2)}}\;dz \right)^{\frac{1}{p}}$ which is finite since $pr'(n-2) < n$. Therefore,
\[
\left\vert \int_{A} F_{k}(y) \;dy \right\vert < c_{4} \left( c_{5} \mu(A)^{\frac{1}{p'}} \right)^{\frac{1}{r'}} \qquad \text{for all}\;\; k \in \mathbb{N}.
\]
Hence, whenever  $\mu(A) < \delta:= \left( \eta^{r'} c_{4}^{-r'} c_{5}^{-1} \right)^{p'} \in \left( 0, \infty \right)$, we have 
$
\left\vert \int_{A} F_{k}(y)\;dy \right\vert < \eta.
$
This proves that the sequence $\lbrace F_{k}\rbrace_{k \in \mathbb{N}}$ is uniformly integrable. Now, applying the Vitali's Convergence Theorem on the ball $B(0,\tilde{R})$, we get 
\[
\int_{B(0,\tilde{R})} F_{k}(y) \;dy \longrightarrow \int_{B(0,\tilde{R})} F_{0}(y) \;dy \qquad \text{as}\;\; k \rightarrow \infty.
\]
That is,
$
I_{k,\tilde{R}} \longrightarrow I_{0,\tilde{R}}$ as $k \rightarrow \infty.$
So there exists $K \in \mathbb{N}$ such that 
$
|I_{k,\tilde{R}} - I_{0,\tilde{R}}|< \frac{\epsilon}{3}$ for $k \geq K.$
Thus, for $k \geq K$ we have
\[
\begin{split}
\left| \int_{\mathbb{R}^n} E_{n}(x_{k}-y)f(y) \;dy - \int_{\mathbb{R}^n} E_{n}(x_{0}-y)f(y) \;dy \right| 
&\leq \left| I_{k,\tilde{R}} - I_{0,\tilde{R}} \right| + \left| II_{{k},\tilde{R}} \right| + \left| II_{0,\tilde{R}} \right| <\epsilon
\end{split} 
\]
Hence,
\[
\int_{\mathbb{R}^n} E_{n}(x_{k}-y)f(y) \;dy \longrightarrow 
\int_{\mathbb{R}^n} E_{n}(x_{0}-y)f(y) \;dy \qquad \text{as} \;\;k\rightarrow \infty.
\]
Since $x_0 \in \mathbb{R}^n$ was arbitrary, this shows that $u$ is continuous on $\mathbb{R}^n$.

\medskip
\medskip
\medskip
\begin{Lem}\label{lemma3}
Let $f: \mathbb{R}^n \rightarrow \mathbb{R}$ be a Lebesgue measurable function. Define 
\[
u(x):= \int_{\mathbb{R}^n} E_{n}(x-y)f(y) \;dy \qquad\qquad \text{for}\;\;x \in \mathbb{R}^n.
\]
If  $\frac{f(y)}{1+|y|^{n-2}} \in L^1(\mathbb{R}^n) $ and $f \in L^{r}_{loc}(\mathbb{R}^n)$ for some $r \in (n,\infty ]$, then $u$ is well-defined in $\mathbb{R}^n$ and belongs to $\mathcal{C}^{1}(\mathbb{R}^n)$. Moreover, for each $j \in \lbrace 1,...,n\rbrace $ and $x \in \mathbb{R}^n$ we have 
 \[
 \partial_{j} u (x) = \int_{\mathbb{R}^n} \partial_{j} E_{n}(x-y)f(y)\;dy.
 \]
\end{Lem}
\noindent{\bf Proof.} By Lemma \ref{lemma2}, $u$ is well-defined and continuous in $\mathbb{R}^n$. Now, we show that $u$ is differentiable in $\mathbb{R}^{n}$.
Fix $j \in \lbrace 1,...,n \rbrace$. Let $x_{0} \in \mathbb{R}^{n}$, and let $ \lbrace h_k \rbrace_{k \in \mathbb{N}} \subset \mathbb{R} \setminus \lbrace 0 \rbrace$ with $h_k \rightarrow 0$ as $k \rightarrow \infty$. Without loss of generality, we may suppose that $\lvert h_k \rvert < 1$ for all $k \in \mathbb{N}$. We will show that
\[
\frac{u(x_{0}+h_{k}e_{j})-u(x_{0})}{h_k} \longrightarrow  \int_{\mathbb{R}^n} \partial_{j} E_{n}(x_{0}-y)f(y)\;dy \qquad \text{as} \;\; k \rightarrow \infty.
\]
For $R > 0$ and $k \in \mathbb{N}$, set
\[
I_{k,R} := \int_{B(0,R)} \frac{E_{n}(x_{0}+h_{k}e_{j}-y) - E_{n}(x_{0}-y)}{h_k}f(y)\;dy, 
\]
\[
II_{k,R} := \int_{\mathbb{R}^n \setminus B(0,R)} \frac{E_{n}(x_{0}+h_{k}e_{j}-y) - E_{n}(x_{0}-y)}{h_k}f(y)\;dy
\]
and 
\[
I_{0,R} := \int_{B(0,R)} (\partial_{j} E_{n})(x_{0}-y)f(y)\;dy, 
\qquad 
II_{0,R} := \int_{\mathbb{R}^n \setminus B(0,R)} (\partial_{j} E_{n})(x_{0}-y)f(y)\;dy.
\]
Note that
\[
\frac{u(x_{0}+h_{k}e_{j})-u(x_{0})}{h_k}  = I_{k,R} + II_{k,R}
\qquad \text{and} \qquad
\int_{\mathbb{R}^n} \partial_{j} E_{n}(x_{0}-y)f(y)\;dy = I_{0,R} + II_{0,R}.
\]
Let $\epsilon > 0$. First, we claim that 
\[
\text{sup}_{k \in \mathbb{N}} \left| II_{k,R} \right| < \frac{\epsilon}{3} \qquad \text{and} \qquad |II_{0,R}| < \frac{\epsilon}{3} \qquad \text{for}\;\;R\;\; \text{large enough}.
\]
For every $R > 2\left( \lvert  x_0 \rvert + 1 \right) > 1,$ $y \in \mathbb{R}^n \setminus B(0,R)$, $k \in \mathbb{N}$ and $t \in (0,1)$, we have 
\[
|x_{0} + th_{k}e_{j} - y|
\geq |y|-\left( |x_{0}| + |th_{k}e_{j}| \right)
 > |y|-\left( |x_{0}| + 1 \right) > \frac{R}{2},
\qquad \text{and} 
\]
\[ 
 |x_{0} + th_{k}e_{j} - y|
  \geq |y|-\left( |x_{0}| + |th_{k}e_{j}| \right) 
   > \frac{|y|}{2} + \frac{R}{2} - \left( |x_{0}| + 1 \right) > \frac{|y|}{2},
\]
which implies that 
\[
|x_{0} + th_{k}e_{j} - y|^{n-1} 
\geq \frac{1}{2} \left[ \left( \frac{R}{2}\right)^{n-1} + \left( \frac{|y|}{2} \right)^{n-1} \right] 
 \geq  \frac{1}{2^{n}} \left(1+|y|^{n-1} \right).
\]
Thus, for $R > 2\left( \lvert  x_0 \rvert + 1 \right)$, applying the Mean Value Theorem, we have 
\[
\begin{split}
\text{sup}_{k \in \mathbb{N}}\lvert II_{k,R} \rvert 
&= \text{sup}_{k \in \mathbb{N}} \left\lvert \int_{\mathbb{R}^n \setminus B(0,R)} \frac{E_{n}(x_{0}+h_{k}e_{j}-y) - E_{n}(x_{0}-y)}{h_k}f(y)\;dy \right\rvert \\
&= \text{sup}_{k \in \mathbb{N}} \left\lvert \int_{\mathbb{R}^n \setminus B(0,R)} \left( \int_{0}^{1} \left( \partial_{j}E_n \right) \left( x_{0}+th_{k}e_{j}-y \right)\;dt \right) f(y)\;dy \right\rvert \\
&\leq \text{sup}_{k \in \mathbb{N}} \int_{\mathbb{R}^n \setminus B(0,R)}  \left( \int_{0}^{1} \left\lvert \left(\partial_{j}E_n \right) \left( x_{0}+th_{k}e_{j}-y \right)\right\rvert \;dt \right)  \left\lvert f(y) \right\rvert\;dy  \\
&\leq \frac{1}{\omega_{n-1}} \text{sup}_{k \in \mathbb{N}}\int_{\mathbb{R}^n \setminus B(0,R)} \left( \int_{0}^{1} \frac{\left\lvert f(y) \right\rvert}{\left\vert x_{0}+th_{k}e_{j}-y \right\vert^{n-1}}  \;dt \right)  \;dy  \\
&\leq \frac{2^{n}}{\omega_{n-1}}\int_{\mathbb{R}^n \setminus B(0,R)} \int_{0}^{1} \frac{\left\lvert f(y) \right\rvert}{1 + |y|^{n-1}}  \;dt\;dy  \\
&= \frac{2^{n}}{\omega_{n-1}}\int_{\mathbb{R}^n \setminus B(0,R)} \frac{\left\lvert f(y) \right\rvert}{1 + |y|^{n-1}} \;dy  \leq \frac{2^{n}}{\omega_{n-1}}\int_{\mathbb{R}^n \setminus B(0,R)} \frac{\left\lvert f(y) \right\rvert}{1 + |y|^{n-2}} \;dy  \\
\end{split}
\]
where the last integral on the right-hand side converges to zero as $R$ goes to infinity, using the Lebesgue Dominated Convergence Theorem. Hence there exists $R_{1} > 0$ such that 
\[
\text{sup}_{k \in \mathbb{N}} |II_{k,R}| < \frac{\epsilon}{3} \qquad \text{for}\;\;R \geq R_1.
\]
Next, we want to claim that $|II_{0,R}| < \frac{\epsilon}{3}$  for $R$ large enough. 

By using the same treat as above, for every $R > \text{max}\left( 2|x_{0}|, 1 \right)$ and $y \in \mathbb{R}^n \setminus B(0,R)$ , we obtain
\[
|x_{0}-y|^{n-1} \geq \frac{1}{2}  \left(1+|y|^{n-1} \right).
\]
Thus, for $R > \text{max}\left( 2|x_{0}|, 1 \right)$, we have
\[
\begin{split}
|II_{0,R}| 
&\leq \int_{\mathbb{R}^n \setminus B(0,R)} \left\vert (\partial_{j} E_{n})(x_{0}-y) \right\vert \left\vert f(y) \right\vert \;dy \\
&\leq \frac{1}{\omega_{n-1}} \int_{\mathbb{R}^n \setminus B(0,R)} \frac{|f(y)|}{|x_{0}-y|^{n-1}} \;dy \\ 
&\leq  \frac{2^{n}}{\omega_{n-1}}\int_{\mathbb{R}^n \setminus B(0,R)} \frac{|f(y)|}{1+ |y|^{n-1}} \;dy \leq  \frac{2^{n}}{\omega_{n-1}}\int_{\mathbb{R}^n \setminus B(0,R)} \frac{|f(y)|}{1+ |y|^{n-2}} \;dy,
\end{split}
\]
where the last integral on the right-hand side converges to zero as $R$ approaches infinity, using the Lebesgue Dominated Convergence Theorem. So, there exists $R_{2} > 0$ such that 
$
|II_{0,R}| < \frac{\epsilon}{3}$ for $R \geq R_2.$
Setting $\tilde{R}:= \text{max}\lbrace R_{1}, R_{2}\rbrace$, then 
\[
\text{sup}_{k \in \mathbb{N}} |II_{k,\tilde{R}}| < \frac{\epsilon}{3} \qquad \text{and} \qquad |II_{0,\tilde{R}}| < \frac{\epsilon}{3}.
\]
Now, we show that 
$
I_{k,\tilde{R}} \longrightarrow I_{0,\tilde{R}}$ as $k \rightarrow \infty,$
using the Vitali's Convergence Theorem on $B(0,\tilde{R})$. For 
$ y \in B(0,\tilde{R})$ and $k \in \mathbb{N}$, set 
\[
F_{k}(y):= \frac{E_{n}(x_{0}+h_{k}e_{j}-y) - E_{n}(x_{0}-y)}{h_k}f(y) \qquad \text{and} \qquad 
F_{0}(y):= (\partial_{j} E_{n})(x_{0}-y)f(y),
\]
which is the kernels of $I_{k,\tilde{R}}$ and $I_{0,\tilde{R}}$, respectively. Clearly that $F_{k} \rightarrow F_{0}$ pointwise almost everywhere on $B(0,\tilde{R})$, as $k \rightarrow \infty$. Now, let us check the uniform integrability of the sequence $\lbrace F_{k}\rbrace_{k \in \mathbb{N}}$. Let $\eta >0$, and let $A \subset B(0,\tilde{R})$. For each $k \in \mathbb{N}$, applying the Mean value theorem, H{\"o}lder's inequality and Minkowski's inequality, we have

\[
\begin{split}
\left\vert \int_{A} F_{k}(y) \;dy \right\vert 
&= \left\lvert \int_{A} \frac{E_{n}(x_{0}+h_{k}e_{j}-y) - E_{n}(x_{0}-y)}{h_k}f(y)\;dy \right\rvert \\
&= \left\lvert \int_{A} \left( \int_{0}^{1} \left( \partial_{j}E_n \right) \left( x_{0}+th_{k}e_{j}-y \right)\;dt \right) f(y)\;dy \right\rvert \\
&\leq \int_{A}  \left( \int_{0}^{1} \left\lvert \left(\partial_{j}E_n \right) \left( x_{0}+th_{k}e_{j}-y \right)\right\rvert \;dt \right)  \left\lvert f(y) \right\rvert\;dy  \\
&\leq \left\Vert f \right\Vert_{L^{r}(B(0,\tilde{R}))} 
      \left\Vert  \int_{0}^{1} \left\lvert \left(\partial_{j}E_n \right) \left( x_{0}+th_{k}e_{j}-y \right)\right\rvert \;dt \right\Vert_{L^{r'}(A)}\\
&\leq \left\Vert f \right\Vert_{L^{r}(B(0,\tilde{R}))}
\int_{0}^{1} \left\Vert \left(\partial_{j}E_n \right) \left( x_{0}+th_{k}e_{j}-y \right)\right\Vert_{L^{r'}(A)} \;dt \\
&\leq \frac{\left\Vert f \right\Vert_{L^{r}(B(0,\tilde{R}))}}{\omega_{n-1}} \int_{0}^{1} \left\Vert  \frac{1}{ \left(x_{0}+th_{k}e_{j}-y \right)^{n-1}}  \right\Vert_{L^{r'}(A)} \;dt \\
&= c_{1} \int_{0}^{1} \left\Vert  \frac{1}{ \left(x_{0}+th_{k}e_{j}-y \right)^{n-1}}  \right\Vert_{L^{r'}(A)} \;dt, \\
\end{split}
\]
where $\frac{1}{r} + \frac{1}{r'} = 1$ and $c_{1}=c_{1}(n,r,\tilde{R}):= \frac{\left\Vert f \right\Vert_{L^{r}(B(0,\tilde{R}))}}{\omega_{n-1}}$ which is finite as $f \in L^{r}_{loc}(\mathbb{R}^n)$. Note that since $r > n$, equivalently, $r'(n-1) < n$, that is, $\frac{n}{r'(n-1)} > 1$. Pick $p \in \left( 1,\frac{n}{r'(n-1)} \right)$. For $k \in \mathbb{N}$ and $t \in (0,1)$, set 
$
x_{t,k} := x_{0}+th_{k}e_{j}.
$
Note that since $t \in (0,1)$ and $h_k \rightarrow 0$ as $k \rightarrow 0$, the sequence $\lbrace x_{t,k} \rbrace$ converges to $x_0$ as $k \rightarrow \infty$ for all $t \in (0,1)$. In particular, there is $M>0$ so that $\left\vert x_{t,k} \right\vert < M $ for all $k \in \mathbb{N}$ and $t \in (0,1)$. Applying H{\"o}lder's inequality again, for any $t \in (0,1)$, we have 
\[
\begin{split}
\left\Vert  \frac{1}{ \left(x_{0}+th_{k}e_{j}-y \right)^{n-1}}  \right\Vert_{L^{r'}(A)}^{r'} 
&= \int_{A} \frac{1}{|x_{t,k}-y|^{r'(n-1)}}\;dy \\
&= \int_{x_{t,k}-A} \frac{1}{|z|^{r'(n-1)}}\;dz \\
&= \int_{B(0,\tilde{R}+M)} \frac{1}{|z|^{r'(n-1)}} \chi_{x_{t,k}-A}(z)\;dz \\
&\leq \left( \int_{B(0,\tilde{R}+M)} \frac{1}{|z|^{pr'(n-1)}}\;dz \right)^{\frac{1}{p}}  \left( \int_{B(0,\tilde{R}+M)} \left( \chi_{x_{t,k}-A}(z)\right)^{p'} \;dz \right)^{\frac{1}{p'}} = c_{2} \left( \mu \left( A \right) \right)^{\frac{1}{p'}} ,
\end{split}
\]
where $\frac{1}{p} + \frac{1}{p'} = 1$ and $c_{2}=c_{2}(n, \tilde{R},M,r,p):= \left( \int_{B(0,\tilde{R}+M)} \frac{1}{|z|^{pr'(n-1)}}\;dz \right)^{\frac{1}{p}}$ which is finite since $pr'(n-1) < n$. Therefore,
\[
\left\vert \int_{A} F_{k}(y) \;dy \right\vert < c_{1}  \int_{0}^{1} \left( c_{2} \mu(A)^{\frac{1}{p'}} \right)^{\frac{1}{r'}}\;dt 
= c_{1} \left( c_{2} \mu(A)^{\frac{1}{p'}} \right)^{\frac{1}{r'}}   \qquad \text{for all}\;\; k \in \mathbb{N}.
\]
Hence, whenever  $\mu(A) < \delta:= \left( \eta^{r'} c_{1}^{-r'} c_{2}^{-1} \right)^{p'} \in \left( 0, \infty \right)$, we have 
$
\left\vert \int_{A} F_{k}(y)\;dy \right\vert < \eta.
$
This proves that the sequence $\lbrace F_{k}\rbrace_{k \in \mathbb{N}}$ is uniformly integrable. Now, applying the Vitali's Convergence Theorem on $B(0,\tilde{R})$, we get 
\[
\int_{B(0,\tilde{R})} F_{k}(y) \;dy \longrightarrow \int_{B(0,\tilde{R})} F_{0}(y) \;dy \qquad \text{as}\;\; k \rightarrow \infty.
\]
That is,
$
I_{k,\tilde{R}} \longrightarrow I_{0,\tilde{R}}$ as $k \rightarrow \infty.
$
So, there exist $K \in \mathbb{N}$ such that 
$
|I_{k,\tilde{R}} - I_{0,\tilde{R}}|< \frac{\epsilon}{3}$ for $k \geq K.$
Thus, for $k \geq K$ we have
\[
\begin{split}
\left\vert \frac{u(x_{0}+h_{k}e_{j})-u(x_{0})}{h_k} - \int_{\mathbb{R}^n} \partial_{j} E_{n}(x_{0}-y)f(y)\;dy \right\vert 
&=  \left| \left( I_{{k},\tilde{R}} + II_{{k},\tilde{R}} \right) - \left( I_{0,\tilde{R}} + II_{0,\tilde{R}} \right) \right| \\
&\leq \left| I_{k,\tilde{R}} - I_{0,\tilde{R}} \right| + \left| II_{{k},\tilde{R}} \right| + \left| II_{0,\tilde{R}} \right| < \epsilon
\end{split} 
\]
Hence,
\[
\frac{u(x_{0}+h_{k}e_{j})-u(x_{0})}{h_k} \longrightarrow  \int_{\mathbb{R}^n} \partial_{j} E_{n}(x_{0}-y)f(y)\;dy \qquad \text{as} \;\; k \rightarrow \infty.
\]
Since $x_0 \in \mathbb{R}^n$ and $j \in \lbrace 1,...,n\rbrace $ were arbitrary, this shows that $u$ is differentiable on $\mathbb{R}^n$. Moreover, for each $j \in \lbrace 1,...,n\rbrace $ and $x \in \mathbb{R}^n$ we have an explicit formula 
 \[
 \partial_{j} u (x) = \int_{\mathbb{R}^n} \partial_{j} E_{n}(x-y)f(y)\;dy.
 \]
Finally, we will show that each $\partial_{j}u$ is continuous on $\mathbb{R}^n$. Fix $j \in \lbrace 1,...,n\rbrace $ and $x_0 \in \mathbb{R}^n$. Let $\lbrace x_k \rbrace_{k \in \mathbb{N}} \subset \mathbb{R}^n$ so that $x_k \rightarrow x_0$ as $k \rightarrow \infty$. We need to show that 
\[
\int_{\mathbb{R}^n} \partial_{j}E_{n}(x_{k}-y)f(y) \;dy \longrightarrow 
\int_{\mathbb{R}^n} \partial_{j}E_{n}(x_{0}-y)f(y) \;dy \qquad \text{as} \;\; 
k\rightarrow \infty.
\]
For $R>0$ and $k \in \mathbb{N}$, set
\[
I_{k,R}:= \int_{B(0,R)} \partial_{j}E_{n}(x_{k}-y)f(y) \;dy, \qquad\;\;
II_{k,R}:=\int_{\mathbb{R}^n \setminus B(0,R)} \partial_{j}E_{n}(x_{k}-y)f(y) \;dy
\]
and 
\[
I_{0,R}:= \int_{B(0,R)} \partial_{j}E_{n}(x_{0}-y)f(y) \;dy, \qquad\;\;
II_{0,R}:=\int_{\mathbb{R}^n \setminus B(0,R)} \partial_{j}E_{n}(x_{0}-y)f(y) \;dy.
\]
Note that 
\[
\int_{\mathbb{R}^n} \partial_{j}E_{n}(x_{k}-y)f(y) \;dy = I_{k,R} + II_{k,R}
\qquad
\text{and} 
\qquad
\int_{\mathbb{R}^n} \partial_{j}E_{n}(x_{0}-y)f(y) \;dy = I_{0,R} + II_{0,R},
\]
also note that since $\lbrace x_{k}\rbrace_{k \in \mathbb{N}}$ is convergent, then it is bounded. We pick $M \in (\frac{1}{2},\infty)$ such that $\lvert x_k \rvert<M$ for all $k \in \mathbb{N}$.
Let $\epsilon > 0$. First, we claim that 
\[
\text{sup}_{k \in \mathbb{N}} \left| II_{k,R} \right| < \frac{\epsilon}{3} \qquad \text{and} \qquad |II_{0,R}| < \frac{\epsilon}{3} \qquad \text{for}\;\;R\;\; \text{large enough}.
\]
For every $R > 2M > 1,$ $y \in \mathbb{R}^n \setminus B(0,R)$ and $k \in \mathbb{N}$ we have \eqref{ineq3} and \eqref{ineq4}
which implies that 
\[
|x_{k}-y|^{n-1} \geq \frac{1}{2} \left[ M^{n-1} + \left( \frac{|y|}{2} \right)^{n-1} \right] 
\geq  \frac{1}{2^{n}} \left(1+|y|^{n-1} \right).
\]
Thus, for $R>2M$ we have
\[
\begin{split}
\text{sup}_{k \in \mathbb{N}} |II_{k,R}|
&\leq \text{sup}_{k \in \mathbb{N}} \int_{\mathbb{R}^n \setminus B(0,R)} \left\vert \partial_{j}E_{n}(x_{k}-y) \right\vert \left\vert f(y) \right\vert \;dy \\ 
&\leq \frac{1}{\omega_{n-1}} \text{sup}_{k \in \mathbb{N}}  \int_{\mathbb{R}^n \setminus B(0,R)} \frac{|f(y)|}{|x_{k}-y|^{n-1}} \;dy \\ 
&\leq  \frac{2^n}{\omega_{n-1}}\int_{\mathbb{R}^n \setminus B(0,R)} \frac{|f(y)|}{1+ |y|^{n-1}} \;dy \leq  \frac{2^n}{\omega_{n-1}}\int_{\mathbb{R}^n \setminus B(0,R)} \frac{|f(y)|}{1+ |y|^{n-2}} \;dy. \\
\end{split}
\]
where the last integral on the right-hand side converges to zero as $R$ goes to infinity, using the Lebesgue Dominated Convergence Theorem. So there exists $R_{1} > 0$ such that 
\[
\text{sup}_{k \in \mathbb{N}} |II_{k,R}| < \frac{\epsilon}{3} \qquad \text{for}\;\;R \geq R_1.
\]
Next, we use the same argument as above to establish that  
$|II_{0,R}| < \frac{\epsilon}{3}$  for $R$ large enough. 
In fact, for every $R >\;\text{max}\;(2|x_{0}|,1)>1$ and $y \in \mathbb{R}^n \setminus B(0,R)$ we have \eqref{ineq1} with $x=x_0$,
which implies that 
\[
|x_{0}-y|^{n-1} \geq \frac{1}{2} \left[ \left( \frac{R}{2} \right)^{n-1} + \left( \frac{|y|}{2} \right)^{n-1} \right] = \frac{1}{2^{n}} \left[ R^{n-1} + |y|^{n-1} \right] \geq  \frac{1}{2^n} \left(1+|y|^{n-1} \right),
\]
Thus, for $R >\;\text{max}\;(2|x_{0}|,1)$ we have
\[
\begin{split}
|II_{0,R}| 
&\leq \int_{\mathbb{R}^n \setminus B(0,R)} \left\vert \partial_{j}E_{n}(x_{0}-y) \right\vert \left\vert f(y) \right\vert \;dy \\ 
&\leq \frac{1}{\omega_{n-1}} \int_{\mathbb{R}^n \setminus B(0,R)} \frac{|f(y)|}{|x_{0}-y|^{n-1}} \;dy \\ 
&\leq  \frac{2^n}{\omega_{n-1}}\int_{\mathbb{R}^n \setminus B(0,R)} \frac{|f(y)|}{1+ |y|^{n-1}} \;dy \\
&\leq  \frac{2^n}{\omega_{n-1}}\int_{\mathbb{R}^n \setminus B(0,R)} \frac{|f(y)|}{1+ |y|^{n-2}} \;dy,
\end{split}
\]
where the last integral on the right-hand side converges to zero as $R$ goes to infinity, using the Lebesgue Dominated Convergence Theorem. Hence there exists $R_{2} > 0$ such that 
\[
|II_{0,R}| < \frac{\epsilon}{3} \qquad \text{for}\;\;R \geq R_2.
\]
Setting $\tilde{R}:= \text{max}\lbrace R_{1}, R_{2}\rbrace$, then 
\[
\text{sup}_{k \in \mathbb{N}} |II_{k,\tilde{R}}| < \frac{\epsilon}{3} \qquad \text{and} \qquad |II_{0,\tilde{R}}| < \frac{\epsilon}{3}.
\]
Finally, we show that 
$
I_{k,\tilde{R}} \longrightarrow I_{0,\tilde{R}} \qquad \text{as}\;\;k \rightarrow \infty,
$
using the Vitali's Convergence Theorem on the ball $B(0,\tilde{R})$. For $ y \in B(0,\tilde{R})$ and $k \in \mathbb{N}$, set
\[
F_{k}(y):= \partial_{j}E_{n}(x_{k}-y)f(y) \qquad \text{and} \qquad F_{0}(y):=  \partial_{j}E_{n}(x_{0}-y)f(y),
\]
which is the kernels of $I_{k,\tilde{R}}$ and $I_{0,\tilde{R}}$ ,respectively Clearly that $F_{k} \rightarrow F_{0}$ pointwise almost everywhere on $B(0,\tilde{R})$, as $k \rightarrow \infty$. Now, let us check the uniform integrability of the sequence $\lbrace F_{k}\rbrace_{k \in \mathbb{N}}$. Let $\eta >0$, and let $A \subset B(0,\tilde{R})$. For each $k \in \mathbb{N}$, applying H{\"o}lder's inequality, we have
\[
\begin{split}
\left\vert \int_{A} F_{k}(y) \;dy \right\vert&\leq \frac{1}{\omega_{n-1}}\int_{A} \frac{|f(y)|}{|x_{k}-y|^{n-1}} \;dy \leq \frac{1}{\omega_{n-1}} \|f\|_{L^{r}(A)} \left\lVert\frac{1}{|x_{k}-y|^{n-1}} \right\rVert_{L^{r'}(A)} \\ 
&\leq \frac{1}{\omega_{n-1}} \|f\|_{L^{r}(B(0,\tilde{R}))} \left\lVert \frac{1}{|x_{k}-y|^{n-1}} \right\rVert_{L^{r'}(A)}\leq c_{3} \left\lVert \frac{1}{|x_{k}-y|^{n-1}} \right\rVert_{L^{r'}(A)},
\end{split}
\]
where $c_{3} = c_{3}(n,r, \tilde{R}) := \frac{1}{\omega_{n-1}} \|f\|_{L^{r}(B(0,\tilde{R}))}$ which is finite as $f \in L^{r}_{loc}(\mathbb{R}^n)$. Note that since $r > n$, equivalently, $r'(n-1) < n$, that is, $\frac{n}{r'(n-1)} > 1$. Pick $p \in \left( 1,\frac{n}{r'(n-1)} \right)$ and $k \in \mathbb{N}$, applying H{\"o}lder's inequality again, we have
\[
\begin{split}
\left\lVert \frac{1}{|x_{k}-y|^{n-1}} \right\rVert_{L^{r'}(A)}^{r'} 
&= \int_{A} \frac{1}{|x_{k}-y|^{r'(n-1)}}\;dy \\
&= \int_{x_{k}-A} \frac{1}{|z|^{r'(n-1)}}\;dz \\
&= \int_{B(0,\tilde{R}+M)} \frac{1}{|z|^{r'(n-1)}} \chi_{x_{k}-A}(z)\;dz \\
&\leq \left( \int_{B(0,\tilde{R}+M)} \frac{1}{|z|^{pr'(n-1)}}\;dz \right)^{\frac{1}{p}}  \left( \int_{B(0,\tilde{R}+M)} \left( \chi_{x_{k}-A}(z)\right)^{p'} \;dz \right)^{\frac{1}{p'}} = c_{4} \left( \mu \left( A \right) \right)^{\frac{1}{p'}} ,
\end{split}
\]
where $\frac{1}{p} + \frac{1}{p'} = 1$ and $c_{4}=c_{4}(n,\tilde{R},M,r,p):= \left( \int_{B(0,\tilde{R}+M)} \frac{1}{|z|^{pr'(n-1)}}\;dz \right)^{\frac{1}{p}}$ which is finite since $pr'(n-1) < n$. Therefore,
\[
\left\vert \int_{A} F_{k}(y) \;dy \right\vert < c_{3} \left( c_{4} \mu(A)^{\frac{1}{p'}} \right)^{\frac{1}{r'}} \qquad \text{for all}\;\; k \in \mathbb{N}.
\]
Hence, whenever  $\mu(A) < \delta:= \left( \eta^{r'} c_{3}^{-r'} c_{4}^{-1} \right)^{p'} \in \left( 0, \infty \right)$, we have 
$
\left\vert \int_{A} F_{k}(y)\;dy \right\vert < \eta.
$
This proves that the sequence $\lbrace F_{k}\rbrace_{k \in \mathbb{N}}$ is uniformly integrable. Now, applying the Vitali's Convergence Theorem on $B(0,\tilde{R})$, we get 
\[
\int_{B(0,\tilde{R})} F_{k}(y) \;dy \longrightarrow \int_{B(0,\tilde{R})} F_{0}(y) \;dy \qquad \text{as}\;\; k \rightarrow \infty.
\]
That is,
$
I_{k,\tilde{R}} \longrightarrow I_{0,\tilde{R}}$ as $k \rightarrow \infty.$
So there exist $K \in \mathbb{N}$ such that 
$
|I_{k,\tilde{R}} - I_{0,\tilde{R}}|< \frac{\epsilon}{3}$ for $k \geq K.$
Thus, for $k \geq K$ we have
\[
\begin{split}
\left| \int_{\mathbb{R}^n} \partial_{j} E_{n}(x_{k}-y)f(y) \;dy - \int_{\mathbb{R}^n} \partial_{j} E_{n}(x_{0}-y)f(y) \;dy \right| 
&=  \left| \left( I_{{k},\tilde{R}} + II_{{k},\tilde{R}} \right) - \left( I_{0,\tilde{R}} + II_{0,\tilde{R}} \right) \right| < \epsilon
\end{split} 
\]
Hence,
\[
\int_{\mathbb{R}^n} \partial_{j} E_{n}(x_{k}-y)f(y) \;dy \longrightarrow 
\int_{\mathbb{R}^n} \partial_{j} E_{n}(x_{0}-y)f(y) \;dy \qquad \text{as} \;\;k\rightarrow \infty.
\]
Since $x_0 \in \mathbb{R}^n$ and $j \in \lbrace 1,...,n \rbrace$ were arbitrary, this shows that each $\partial_{j}u$ is continuous on $\mathbb{R}^n$.
\\
\\
Now, we are ready to prove the main theorem.

\medskip
\noindent{\bf Proof of Theorem \ref{mainthm}}
In view of Lemma \ref{lemma1}, Lemma \ref{lemma2} and Lemma \ref{lemma3}, $u$ is well-defined in $\mathbb{R}^n$, and is  of class $\mathcal{C}^{1}(\mathbb{R}^n)$. Moreover, for each $j \in \lbrace 1,...,n\rbrace $ and $x \in \mathbb{R}^n$ we have 
 \[
 \partial_{j} u (x) = \int_{\mathbb{R}^n} \partial_{j} E_{n}(x-y)f(y)\;dy = \int_{\mathbb{R}^n} \partial_{j} E_{n}(y)f(x-y)\;dy.
 \] Now, we show that $u$ is twice differentiable in $\mathbb{R}^{n}$. Fix $i,j \in \lbrace 1,...,n \rbrace$. Let $x_{0} \in \mathbb{R}^{n}$, and let $ \lbrace h_k \rbrace_{k \in \mathbb{N}} \subset \mathbb{R} \setminus \lbrace 0 \rbrace$ with $h_k \rightarrow 0$ as $k \rightarrow \infty$. Without loss of generality, we may suppose that $\lvert h_k \rvert < 1$ for all $k \in \mathbb{N}$. We will show that
\[
\frac{ \partial_{j}u(x_{0}+h_{k}e_{i})-\partial_{j}u(x_{0})}{h_k} \longrightarrow  \int_{\mathbb{R}^n} \partial_{j} E_{n}(x_{0}-y)\partial_{i}f(y)\;dy \qquad \text{as} \;\; k \rightarrow \infty.
\]
For $R > 0$ and $k \in \mathbb{N}$, setting
\[
I_{k,R} := \int_{B(0,R)} \int_{0}^{1} \partial_{j}E_{n}\left( x_{0}+th_{k}e_{i}-y \right) \partial_{i}f\left( y \right)\;dt\;dy,
\] 
\[
II_{k,R} := \int_{\mathbb{R}^n \setminus B(0,R)}  \int_{0}^{1} \partial_{j}E_{n}\left( x_{0}+th_{k}e_{i}-y \right) \partial_{i}f\left( y \right)\;dt\;dy,
\]
and 
\[
I_{0,R} := \int_{B(0,R)} \partial_{j} E_{n}(x_{0}-y)\partial_{i}f(y)\;dy, 
\qquad 
II_{0,R} := \int_{\mathbb{R}^n \setminus B(0,R)} \partial_{j} E_{n}(x_{0}-y)\partial_{i}f(y)\;dy.
\]
By the Mean Value Theorem and change of variables, we have 
\[
\frac{\partial_{i}u(x_{0}+h_{k}e_{i})-\partial_{i}u(x_{0})}{h_k}  = I_{k,R} + II_{k,R}
\qquad \text{and} \qquad
\int_{\mathbb{R}^n} \partial_{j} E_{n}(x_{0}-y)\partial_{i}f(y)\;dy = I_{0,R} + II_{0,R}.
\]
Let $\epsilon > 0$. First, we claim that 
\[
\text{sup}_{k \in \mathbb{N}} \left| II_{k,R} \right| < \frac{\epsilon}{3} \qquad \text{and} \qquad |II_{0,R}| < \frac{\epsilon}{3} \qquad \text{for}\;\;R\;\; \text{large enough}.
\]
For every $R > 2\left( \lvert  x_0 \rvert + 1 \right) > 1,$ $y \in \mathbb{R}^n \setminus B(0,R)$, $k \in \mathbb{N}$ and $t \in (0,1)$, we have 
\[
|x_{0} + th_{k}e_{i} - y| \geq |y|-|x_{0} + th_{k}e_{i}|  
\geq |y|-\left( |x_{0}| + |th_{k}e_{i}| \right) > |y|-\left( |x_{0}| + 1 \right) > R - \frac{R}{2} = \frac{R}{2}, 
\qquad \text{and} 
\]
\[ 
 |x_{0} + th_{k}e_{i} - y| \geq |y|-|x_{0} + th_{k}e_{i}| \geq |y|-\left( |x_{0}| + |th_{k}e_{i}| \right) = \frac{|y|}{2} + \frac{|y|}{2} - \left( |x_{0}| + |th_{k}e_{i}| \right) > \frac{|y|}{2} + \frac{R}{2} - \left( |x_{0}| + 1 \right) > \frac{|y|}{2},
\]
which implies that 
\[
|x_{0} + th_{k}e_{i} - y|^{n-1} 
\geq \frac{1}{2} \left[ \left( \frac{R}{2}\right)^{n-1} + \left( \frac{|y|}{2} \right)^{n-1} \right] 
= \frac{1}{2^{n}} \left[ R^{n-1} + |y|^{n-1} \right] \geq  \frac{1}{2^{n}} \left(1+|y|^{n-1} \right),
\]
Thus, for $R > 2\left( \lvert  x_0 \rvert + 1 \right)$, we have 
\[
\begin{split}
\text{sup}_{k \in \mathbb{N}}\lvert II_{k,R} \rvert 
&= \text{sup}_{k \in \mathbb{N}} \left\lvert \int_{\mathbb{R}^n \setminus B(0,R)}  \int_{0}^{1} \partial_{j}E_n  \left( x_{0}+th_{k}e_{i}-y \right) \partial_{i}f(y)\;dt\;dy \right\rvert \\
&\leq \text{sup}_{k \in \mathbb{N}} \int_{\mathbb{R}^n \setminus B(0,R)}  \left( \int_{0}^{1} \left\lvert \partial_{j}E_n \left( x_{0}+th_{k}e_{i}-y \right)\right\rvert \;dt \right)  \left\lvert \partial_{i}f(y) \right\rvert\;dy  \\
&\leq \frac{1}{\omega_{n-1}} \text{sup}_{k \in \mathbb{N}}\int_{\mathbb{R}^n \setminus B(0,R)} \left( \int_{0}^{1} \frac{\left\lvert \partial_{i}f(y) \right\rvert}{\left\vert x_{0}+th_{k}e_{i}-y \right\vert^{n-1}}  \;dt \right)  \;dy  \\
&\leq \frac{2^{n}}{\omega_{n-1}}\int_{\mathbb{R}^n \setminus B(0,R)} \int_{0}^{1} \frac{\left\lvert \partial_{i}f(y) \right\rvert}{1 + |y|^{n-1}}  \;dt\;dy  \\
&= \frac{2^{n}}{\omega_{n-1}}\int_{\mathbb{R}^n \setminus B(0,R)} \frac{\left\lvert \partial_{i}f(y) \right\rvert}{1 + |y|^{n-1}} \;dy \leq \frac{2^{n}}{\omega_{n-1}}\int_{\mathbb{R}^n \setminus B(0,R)} \frac{\left\lvert \nabla f(y) \right\rvert}{1 + |y|^{n-1}} \;dy  \\
\end{split}
\]
where the last integral on the right-hand side converges to zero as $R$ goes to infinity, using the Lebesgue Dominated Convergence Theorem. Hence there exists $R_{1} > 0$ such that 
\[
\text{sup}_{k \in \mathbb{N}} |II_{k,R}| < \frac{\epsilon}{3} \qquad \text{for}\;\;R \geq R_1.
\]
Next, we use the same argument as above to establish that  
$|II_{0,R}| < \frac{\epsilon}{3}$  for  $R$ large enough. 
In fact, for every $R > \text{max}\left( 2|x_{0}|, 1 \right)$ and $y \in \mathbb{R}^n \setminus B(0,R)$ we have \eqref{ineq1} with $x=x_0$,
which implies that 
\[
|x_{0}-y|^{n-1} \geq \frac{1}{2} \left[ \left(  \frac{R}{2} \right)^{n-1} + \left( \frac{|y|}{2} \right)^{n-1} \right] = \frac{1}{2^{n}} \left[ R^{n-1} + |y|^{n-1} \right] \geq  \frac{1}{2^{n}} \left(1+|y|^{n-1} \right),
\]
Thus, for $R > \text{max}\left( 2|x_{0}|, 1 \right)$, we have
\[
\begin{split}
|II_{0,R}| 
&\leq \int_{\mathbb{R}^n \setminus B(0,R)} \left\vert (\partial_{j} E_{n})(x_{0}-y) \right\vert \left\vert \partial_{i}f(y) \right\vert \;dy \\
&\leq \frac{1}{\omega_{n-1}} \int_{\mathbb{R}^n \setminus B(0,R)} \frac{|\partial_{i}f(y)|}{|x_{0}-y|^{n-1}} \;dy \\ 
&\leq  \frac{2^{n}}{\omega_{n-1}}\int_{\mathbb{R}^n \setminus B(0,R)} \frac{|\partial_{i}f(y)|}{1+ |y|^{n-1}} \;dy\leq  \frac{2^{n}}{\omega_{n-1}}\int_{\mathbb{R}^n \setminus B(0,R)} \frac{|\nabla f(y)|}{1+ |y|^{n-1}} \;dy,
\end{split}
\]
where the last integral on the right-hand side converges to zero as $R$ goes to infinity, using the Lebesgue Dominated Convergence Theorem. So, there exists $R_{2} > 0$ such that 
$
|II_{0,R}| < \frac{\epsilon}{3}$ for $R \geq R_2.$
Setting $\tilde{R}:= \text{max}\lbrace R_{1}, R_{2}\rbrace$, then 
$
\sup\limits_{k \in \mathbb{N}} |II_{k,\tilde{R}}| < \frac{\epsilon}{3}$ and $|II_{0,\tilde{R}}| < \frac{\epsilon}{3}.$
Now, we show that 
$
I_{k,\tilde{R}} \longrightarrow I_{0,\tilde{R}}$ as $k \rightarrow \infty,$
using the Vitali's Convergence Theorem on the ball $B(0,\tilde{R})$. For 
$ y \in B(0,\tilde{R})$ and $k \in \mathbb{N}$, set 
\[
F_{k}(y):= \int_{0}^{1} \partial_{j}E_{n}\left( x_{0}+th_{k}e_{i}-y \right) \partial_{i}f\left( y \right)\;dt \qquad \text{and} \qquad 
F_{0}(y):= (\partial_{j} E_{n})(x_{0}-y)\partial_{i}f(y),
\]
which is the kernels of $I_{k,\tilde{R}}$ and $I_{0,\tilde{R}}$, respectively. Since $\partial_{j}E_{n}$ is continuous and $[0,1]$ is compact, then $\partial_{j}E_n(x_0 + \cdot h_k e_i - y) \in L^{1}([0,1])$. By the Lebesgue Dominated Convergence Theorem, we have $F_{k} \rightarrow F_{0}$ pointwise almost everywhere on $B(0,\tilde{R})$, as $k \rightarrow \infty$. Now, let us check the uniform integrability of the sequence $\lbrace F_{k}\rbrace_{k \in \mathbb{N}}$. Let $\eta >0$, and let $A \subset B(0,\tilde{R})$. For each $k \in \mathbb{N}$, we have
\[
\begin{split}
\left\vert \int_{A} F_{k}(y) \;dy \right\vert 
&= \left\lvert \int_{A}  \int_{0}^{1} \left( \partial_{j}E_n \right) \left( x_{0}+th_{k}e_{i}-y \right)\partial_{i}f(y)\;dt\;dy \right\rvert \\
&\leq \int_{A} \int_{0}^{1} \left\lvert \left(\partial_{j}E_n \right) \left( x_{0}+th_{k}e_{i}-y \right)\right\rvert \;dt  \left\lvert \partial_{i}f(y) \right\rvert\;dy  \\
&\leq \left\Vert \partial_{i}f \right\Vert_{L^{\infty}(B(0,\tilde{R}))} 
  \int_{0}^{1} \int_{A}\left\lvert \left(\partial_{j}E_n \right) \left( x_{0}+th_{k}e_{i}-y \right)\right\rvert \;dy\;dt \\
&= \left\Vert \partial_{i}f \right\Vert_{L^{\infty}(B(0,\tilde{R}))}
\int_{0}^{1} \left\Vert \left(\partial_{j}E_n \right) \left( x_{0}+th_{k}e_{i}-y \right)\right\Vert_{L^{1}(A)} \;dt \\
&\leq \frac{\left\Vert \partial_{i}f \right\Vert_{L^{\infty}(B(0,\tilde{R}))}}{\omega_{n-1}} \int_{0}^{1} \left\Vert  \frac{1}{ \left(x_{0}+th_{k}e_{i}-y \right)^{n-1}}  \right\Vert_{L^{1}(A)} \;dt \\
&= c_{1} \int_{0}^{1} \left\Vert  \frac{1}{ \left(x_{0}+th_{k}e_{i}-y \right)^{n-1}}  \right\Vert_{L^{1}(A)} \;dt, \\
\end{split}
\]
where $c_{1}=c_{1}(n, \tilde{R}):= \frac{\left\Vert \partial_{i}f \right\Vert_{L^{\infty}(B(0,\tilde{R}))}}{\omega_{n-1}}$ which is finite as $f \in \mathcal{C}^{1}(\mathbb{R}^n)$. Pick $p \in \left( 1,\frac{n}{n-1} \right)$. For $k \in \mathbb{N}$ and $t \in (0,1)$, set 
$
x_{t,k} := x_{0}+th_{k}e_{i}.
$
Note that since $t \in (0,1)$ and $h_k \rightarrow 0$ as $k \rightarrow 0$, the sequence $\lbrace x_{t,k} \rbrace$ converges to $x_0$ as $k \rightarrow \infty$ for all $t \in (0,1)$. In particular, there is $M>0$ so that $\left\vert x_{t,k} \right\vert < M $ for all $k \in \mathbb{N}$ and $t \in (0,1)$. Applying H{\"o}lder's inequality again, for any $t \in (0,1)$ we have 
\[
\begin{split}
\left\Vert  \frac{1}{ \left(x_{0}+th_{k}e_{i}-y \right)^{n-1}}  \right\Vert_{L^{1}(A)} 
&= \int_{A} \frac{1}{|x_{t,k}-y|^{n-1}}\;dy \\
&= \int_{x_{t,k}-A} \frac{1}{|z|^{n-1}}\;dz \\
&= \int_{B(0,\tilde{R}+M)} \frac{1}{|z|^{n-1}} \chi_{x_{t,k}-A}(z)\;dz \\
&\leq \left( \int_{B(0,\tilde{R}+M)} \frac{1}{|z|^{p(n-1)}}\;dz \right)^{\frac{1}{p}}  \left( \int_{B(0,\tilde{R}+M)} \left( \chi_{x_{t,k}-A}(z)\right)^{p'} \;dz \right)^{\frac{1}{p'}} \\
&= c_{2} \left( \mu \left( A \right) \right)^{\frac{1}{p'}} ,
\end{split}
\]
where $\frac{1}{p} + \frac{1}{p'} = 1$ and $c_{2}=c_{2}(n, \tilde{R},M,p):= \left( \int_{B(0,\tilde{R}+M)} \frac{1}{|z|^{p(n-1)}}\;dz \right)^{\frac{1}{p}}$ which is finite since $p(n-1) < n$. Therefore,
\[
\left\vert \int_{A} F_{k}(y) \;dy \right\vert < c_{1}  \int_{0}^{1}  c_{2} \mu(A)^{\frac{1}{p'}} \;dt 
= c_{1}  c_{2} \mu(A)^{\frac{1}{p'}}   \qquad \text{for all}\;\; k \in \mathbb{N}.
\]
Hence, whenever  $\mu(A) < \delta:= \left( \eta c_{1}^{-1} c_{2}^{-1} \right)^{p'} \in \left( 0, \infty \right)$, we have 
$
\left\vert \int_{A} F_{k}(y)\;dy \right\vert < \eta.
$
This proves that the sequence $\lbrace F_{k}\rbrace_{k \in \mathbb{N}}$ is uniformly integrable. Now, applying the Vitali's Convergence Theorem on $B(0,\tilde{R})$, we get 
\[
\int_{B(0,\tilde{R})} F_{k}(y) \;dy \longrightarrow \int_{B(0,\tilde{R})} F_{0}(y) \;dy \qquad \text{as}\;\; k \rightarrow \infty.
\]
That is,
$
I_{k,\tilde{R}} \longrightarrow I_{0,\tilde{R}} \qquad \text{as}\;\;k \rightarrow \infty.
$
So, there exist $K \in \mathbb{N}$ such that 
$
|I_{k,\tilde{R}} - I_{0,\tilde{R}}|< \frac{\epsilon}{3}$ for $ k \geq K.$
Thus, for $k \geq K$ we have
\[
\begin{split}
\left\vert \frac{ \partial_{j}u(x_{0}+h_{k}e_{i})-\partial_{j}u(x_{0})}{h_k} - \int_{\mathbb{R}^n} \partial_{j} E_{n}(x_{0}-y)\partial_{i}f(y)\;dy \right\vert 
&=  \left| \left( I_{{k},\tilde{R}} + II_{{k},\tilde{R}} \right) - \left( I_{0,\tilde{R}} + II_{0,\tilde{R}} \right) \right| <\epsilon.
\end{split} 
\]
Hence,
\[
\frac{\partial_{j}u(x_{0}+h_{k}e_{j})-\partial_{j}u(x_{0})}{h_k} \longrightarrow  \int_{\mathbb{R}^n} \partial_{j} E_{n}(x_{0}-y)\partial_{i}f(y)\;dy \qquad \text{as} \;\; k \rightarrow \infty.
\]
Since $x_0 \in \mathbb{R}^n$ and $i,j \in \lbrace 1,...,n\rbrace $ were arbitrary, this shows that $u$ is twice differentiable in $\mathbb{R}^n$. Moreover, for each $i,j \in \lbrace 1,...,n\rbrace $ and $x \in \mathbb{R}^n$ we have an explicit formula 
 \[
 \partial_{i}\partial_{j} u (x) = \int_{\mathbb{R}^n} \partial_{j} E_{n}(x-y)\partial_{i}f(y)\;dy.
 \]
Now, we will show that each $\partial_{j}\partial_{i}u$ is continuous on $\mathbb{R}^n$. Fix $i,j \in \lbrace 1,...,n\rbrace $ and $x_0 \in \mathbb{R}^n$. Let $\lbrace x_k \rbrace_{k \in \mathbb{N}} \subset \mathbb{R}^n$ so that $x_k \rightarrow x_0$ as $k \rightarrow \infty$. We need to show that 
\[
\int_{\mathbb{R}^n} \partial_{j}E_{n}(x_{k}-y)\partial_{i}f(y) \;dy \longrightarrow 
\int_{\mathbb{R}^n} \partial_{j}E_{n}(x_{0}-y)\partial_{i}f(y) \;dy \qquad \text{as} \;\; 
k\rightarrow \infty.
\]
For $R>0$ and $k \in \mathbb{N}$, set
\[
I_{k,R}:= \int_{B(0,R)} \partial_{j}E_{n}(x_{k}-y)\partial_{i}f(y) \;dy, \qquad\;\;
II_{k,R}:=\int_{\mathbb{R}^n \setminus B(0,R)} \partial_{j}E_{n}(x_{k}-y)\partial_{i}f(y) \;dy
\]
and 
\[
I_{0,R}:= \int_{B(0,R)} \partial_{j}E_{n}(x_{0}-y)\partial_{i}f(y) \;dy, \qquad\;\;
II_{0,R}:=\int_{\mathbb{R}^n \setminus B(0,R)} \partial_{j}E_{n}(x_{0}-y)\partial_{i}f(y) \;dy.
\]
Note that 
\[
\int_{\mathbb{R}^n} \partial_{j}E_{n}(x_{k}-y)\partial_{i}f(y) \;dy = I_{k,R} + II_{k,R}
\qquad
\text{and} 
\qquad
\int_{\mathbb{R}^n} \partial_{j}E_{n}(x_{0}-y)\partial_{i}f(y) \;dy = I_{0,R} + II_{0,R},
\]
also note that since $\lbrace x_{k}\rbrace_{k \in \mathbb{N}}$ is convergent, then it is bounded. We pick $M \in (\frac{1}{2},\infty)$ such that $\lvert x_k \rvert<M$ for all $k \in \mathbb{N}$.
Let $\epsilon > 0$. First, we claim that 
\[
\text{sup}_{k \in \mathbb{N}} \left| II_{k,R} \right| < \frac{\epsilon}{3} \qquad \text{and} \qquad |II_{0,R}| < \frac{\epsilon}{3} \qquad \text{for}\;\;R\;\; \text{large enough}.
\]
For every $R > 2M > 1,$ $y \in \mathbb{R}^n \setminus B(0,R)$ and $k \in \mathbb{N}$ we have \eqref{ineq3} and \eqref{ineq4} 
which implies that 
\[
|x_{k}-y|^{n-1} \geq \frac{1}{2} \left[ M^{n-1} + \left( \frac{|y|}{2} \right)^{n-1} \right] = \frac{1}{2^{n}} \left[ (2M)^{n-1} + |y|^{n-1} \right] \geq  \frac{1}{2^{n}} \left(1+|y|^{n-1} \right).
\]
Thus, for $R>2M$ we have
\[
\begin{split}
\text{sup}_{k \in \mathbb{N}} |II_{k,R}|
&\leq \text{sup}_{k \in \mathbb{N}} \int_{\mathbb{R}^n \setminus B(0,R)} \left\vert \partial_{j}E_{n}(x_{k}-y) \right\vert \left\vert \partial_{i}f(y) \right\vert \;dy \\ 
&\leq \frac{1}{\omega_{n-1}} \text{sup}_{k \in \mathbb{N}}  \int_{\mathbb{R}^n \setminus B(0,R)} \frac{|\partial_{i}f(y)|}{|x_{k}-y|^{n-1}} \;dy \\ 
&\leq  \frac{2^n}{\omega_{n-1}}\int_{\mathbb{R}^n \setminus B(0,R)} \frac{|\partial_{i}f(y)|}{1+ |y|^{n-1}} \;dy \leq  \frac{2^n}{\omega_{n-1}}\int_{\mathbb{R}^n \setminus B(0,R)} \frac{|\nabla f(y)|}{1+ |y|^{n-1}} \;dy. \\
\end{split}
\]
where the last integral on the right-hand side converges to zero as $R$ goes to infinity, using the Lebesgue Dominated Convergence Theorem. So there exists $R_{1} > 0$ such that 
\[
\text{sup}_{k \in \mathbb{N}} |II_{k,R}| < \frac{\epsilon}{3} \qquad \text{for}\;\;R \geq R_1.
\]
Next, we use the same argument as above to establish that  
$
|II_{0,R}| < \frac{\epsilon}{3} $ for $R$ large enough.
In fact, for every $R >\;\text{max}\;(2|x_{0}|,1)>1$ and $y \in \mathbb{R}^n \setminus B(0,R)$ we have \eqref{ineq1} and \eqref{ineq2} with $x=x_0$,
which imply that 
\[
|x_{0}-y|^{n-1} \geq \frac{1}{2} \left[ \left( \frac{R}{2} \right)^{n-1} + \left( \frac{|y|}{2} \right)^{n-1} \right] = \frac{1}{2^{n}} \left[ R^{n-1} + |y|^{n-1} \right] \geq  \frac{1}{2^n} \left(1+|y|^{n-1} \right),
\]
Thus, for $R >\;\text{max}\;(2|x_{0}|,1)$ we have
\[
\begin{split}
|II_{0,R}| 
&\leq \int_{\mathbb{R}^n \setminus B(0,R)} \left\vert \partial_{j}E_{n}(x_{0}-y) \right\vert \left\vert \partial_{i}f(y) \right\vert \;dy \\ 
&\leq \frac{1}{\omega_{n-1}} \int_{\mathbb{R}^n \setminus B(0,R)} \frac{|\partial_{i}f(y)|}{|x_{0}-y|^{n-1}} \;dy \\ 
&\leq  \frac{2^n}{\omega_{n-1}}\int_{\mathbb{R}^n \setminus B(0,R)} \frac{|\partial_{i}f(y)|}{1+ |y|^{n-1}} \;dy \leq  \frac{2^n}{\omega_{n-1}}\int_{\mathbb{R}^n \setminus B(0,R)} \frac{|\nabla f(y)|}{1+ |y|^{n-1}} \;dy,
\end{split}
\]
where the last integral on the right-hand side converges to zero as $R$ goes to infinity, using the Lebesgue Dominated Convergence Theorem. Hence there exists $R_{2} > 0$ such that 
\[
|II_{0,R}| < \frac{\epsilon}{3} \qquad \text{for}\;\;R \geq R_2.
\]
Setting $\tilde{R}:= \text{max}\lbrace R_{1}, R_{2}\rbrace$, then 
\[
\text{sup}_{k \in \mathbb{N}} |II_{k,\tilde{R}}| < \frac{\epsilon}{3} \qquad \text{and} \qquad |II_{0,\tilde{R}}| < \frac{\epsilon}{3}.
\]
Finally, since $\partial_{j}E_{n}$ and $\partial_{i}f$ are continuous and $B(0,\tilde{R})$ is compact, then $\partial_{j}E_{n}(x_k - \cdot)\partial_{i}f(\cdot) \in L^{1}( B(0,\tilde{R}) )$. By the Lebesgue Dominated Convergence Theorem, we have 
$
I_{k,\tilde{R}} \longrightarrow I_{0,\tilde{R}} \; \text{as}\;\;k \rightarrow \infty.
$
(Note that in this case, we don't need to invoke Vitali's Convergence Theorem as in the previous part since $f \in \mathcal{C}^{1}(\mathbb{R}^n)$.) So there exist $K \in \mathbb{N}$ such that 
$
|I_{k,\tilde{R}} - I_{0,\tilde{R}}|< \frac{\epsilon}{3}$ for  $k \geq K$.
Thus, for $k \geq K$ we have
\[
\left| \int_{\mathbb{R}^n} \partial_{j} E_{n}(x_{k}-y)\partial_{i}f(y) \;dy - \int_{\mathbb{R}^n} \partial_{j} E_{n}(x_{0}-y)\partial_{i}f(y) \;dy \right| 
=  \left| \left( I_{{k},\tilde{R}} + II_{{k},\tilde{R}} \right) - \left( I_{0,\tilde{R}} + II_{0,\tilde{R}} \right) \right| \\
<\epsilon.
\]
Hence,
\[
\int_{\mathbb{R}^n} \partial_{j} E_{n}(x_{k}-y)\partial_{i}f(y) \;dy \longrightarrow 
\int_{\mathbb{R}^n} \partial_{j} E_{n}(x_{0}-y)\partial_{i}f(y) \;dy \qquad \text{as} \;\;k\rightarrow \infty.
\]
Since $x_0 \in \mathbb{R}^n$ and $i,j \in \lbrace 1,...,n \rbrace$ were arbitrary, this shows that each $\partial_{i}\partial_{j}u$ is continuous on $\mathbb{R}^n$. Lastly, we will show that the function $u$ solves the Poisson equation $\Delta u = f $ in  $\mathbb{R}^n$.
Let $x_0 \in \mathbb{R}^n$. By integration by parts, we have 
\[
\begin{split}
\Delta u (x_0)  = \sum_{j=1}^{n} \partial_{j}^{2}u(x_0) 
& = \sum_{j=1}^{n} \int_{\mathbb{R}^n} \partial_{j}E_n(x_0 - y) \partial_{j}f(y)\;dy \\
& = \sum_{j=1}^{n} \lim_{\substack{R\rightarrow\infty\\\epsilon \rightarrow 0^+}} \int_{B(x_0,R) \setminus \overline{B(x_0, \epsilon)}} \partial_{j}E_n(x_0 - y) \partial_{j}f(y)\;dy \\
& = \sum_{j=1}^{n}  \lim_{\substack{R\rightarrow\infty\\\epsilon \rightarrow 0^+}} \int_{B(x_0,R) \setminus \overline{B(x_0, \epsilon)}} - \partial_{j}^{2}E_n(x_0 - y) f(y)\;dy\;+  \\
&\quad\sum_{j=1}^{n}\lim_{\substack{R\rightarrow\infty\\\epsilon \rightarrow 0^+}} \int_{ \partial (B(x_0,R) \setminus \overline{B(x_0, \epsilon)})} \partial_{j}E_n(x_0 - y) f(y)\nu_{j}(y) \;d\sigma(y) , 
\end{split}
\]
Since $E_n \in \mathcal{C}^{\infty}(\mathbb{R}^n \setminus \lbrace 0 \rbrace)$ with $\Delta E_n = 0$ on $\mathbb{R}^n \setminus \lbrace 0 \rbrace$ and $\partial_{j}E_n(x)= \frac{1}{\omega_{n-1}} \frac{x_j}{\left\vert x \right\vert^n}$ for all $x \in \mathbb{R}^n \setminus \lbrace 0 \rbrace$, then 
\[
\begin{split}
\Delta u (x_0) 
&= \lim_{\substack{R\rightarrow\infty\\\epsilon \rightarrow 0^+}} \int_{ \partial (B(x_0,R) \setminus \overline{B(x_0, \epsilon)})} 
\sum_{j=1}^{n}\partial_{j}E_n(x_0 - y) f(y)\nu_{j}(y) \;d\sigma(y) \\
&= \lim_{R\rightarrow\infty} \int_{ \partial B(x_0,R)} \sum_{j=1}^{n}\partial_{j}E_n(x_0 - y) f(y)\left(\frac{y-x_0}{R}\right)_{j} \;d\sigma_{R}(y)\;- \\
&\quad \lim_{\epsilon\rightarrow 0^+} \int_{ \partial B(x_0,\epsilon)} \sum_{j=1}^{n}\partial_{j}E_n(x_0 - y) f(y)\left( \frac{y-x_0}{\epsilon}\right)_{j} \;d\sigma_{\epsilon}(y) \\
&= I + II,
\end{split} 
\]
where 
\[
I := \lim_{R\rightarrow\infty} \int_{ \partial B(x_0,R)} \sum_{j=1}^{n}\partial_{j}E_n(x_0 - y) f(y)\left(\frac{y-x_0}{R}\right)_{j} \;d\sigma_{R}(y),
\]
and
\[
II := -\lim_{\epsilon\rightarrow 0^+} \int_{ \partial B(x_0,\epsilon)} \sum_{j=1}^{n}\partial_{j}E_n(x_0 - y) f(y)\left( \frac{y-x_0}{\epsilon}\right)_{j} \;d\sigma_{\epsilon}(y).
\]
Now, we consider the second term $II$. 
\[
\begin{split}
II 
&= -\lim_{\epsilon\rightarrow 0^+} \int_{ \partial B(x_0,\epsilon)} \sum_{j=1}^{n}\partial_{j}E_n(x_0 - y) f(y)\left( \frac{y-x_0}{\epsilon}\right)_{j} \;d\sigma_{\epsilon}(y) \\
&= -\frac{1}{\omega_{n-1}}\lim_{\epsilon\rightarrow 0^+} \int_{ \partial B(x_0,\epsilon)} \sum_{j=1}^{n} \frac{\left(x_0 - y\right)_{j}}{\lvert x_0 - y \rvert^{n}} f(y)\left( \frac{y-x_0}{\epsilon}\right)_{j} \;d\sigma_{\epsilon}(y) \\
&= -\frac{1}{\omega_{n-1}}\lim_{\epsilon\rightarrow 0^+} \frac{1}{\epsilon^{n+1}} \int_{ \partial B(x_0,\epsilon)} \sum_{j=1}^{n} \left(x_0 - y\right)_{j}  \left(y-x_0 \right)_{j} f(y)\;d\sigma_{\epsilon}(y) \\
&= \lim_{\epsilon\rightarrow 0^+} \frac{1}{\omega_{n-1}\epsilon^{n-1}} \int_{ \partial B(x_0,\epsilon)} f(y)\;d\sigma_{\epsilon}(y) \\
&= \lim_{\epsilon\rightarrow 0^+} \frac{1}{\omega_{n-1}\epsilon^{n-1}} \int_{ \partial B(x_0,\epsilon)} f(y) - f(x_0) + f(x_0)\;d\sigma_{\epsilon}(y) \\
&= \left( \lim_{\epsilon\rightarrow 0^+} \frac{1}{\omega_{n-1}\epsilon^{n-1}} \int_{ \partial B(x_0,\epsilon)} f(y) - f(x_0)\;d\sigma_{\epsilon}(y) \right)
+ \left(\lim_{\epsilon\rightarrow 0^+} \frac{1}{\omega_{n-1}\epsilon^{n-1}} \int_{ \partial B(x_0,\epsilon)} f(x_0)\;d\sigma_{\epsilon}(y) \right) \\
&= \left( \lim_{\epsilon\rightarrow 0^+} \frac{1}{\omega_{n-1}\epsilon^{n-1}} \int_{ \partial B(x_0,\epsilon)} f(y) - f(x_0)\;d\sigma_{\epsilon}(y) \right) + f(x_0)= II' + f(x_0),
\end{split}
\]
where 
$
II' := \lim_{\epsilon\rightarrow 0^+} \frac{1}{\omega_{n-1}\epsilon^{n-1}} \int_{ \partial B(x_0,\epsilon)} f(y) - f(x_0)\;d\sigma_{\epsilon}(y).
$
We see that 
\[
\lvert II' \rvert 
\leq \limsup_{\substack{\epsilon \rightarrow 0^+\\ \epsilon \leq 1}}\frac{1}{\omega_{n-1}\epsilon^{n-1}} \int_{ \partial B(x_0,\epsilon)} \epsilon  \sup_{\overline{B(x_0,\epsilon)}} \lvert \nabla f \rvert \;d\sigma_{\epsilon}(y) 
= \limsup_{\substack{\epsilon \rightarrow 0^+\\ \epsilon \leq 1}}
\left( \epsilon  \sup_{\overline{B(x_0,\epsilon)}} \lvert \nabla f \rvert \right) 
= 0.
\]
Hence, $II = f(x_0)$. Next, we consider the first term $I$.
\[
\begin{split}
I 
&= \lim_{R\rightarrow\infty} \int_{ \partial B(x_0,R)} \sum_{j=1}^{n}\partial_{j}E_n(x_0 - y) f(y)\left(\frac{y-x_0}{R}\right)_{j} \;d\sigma_{R}(y) \\
&= \lim_{R\rightarrow \infty} \frac{1}{\omega_{n-1}R^{n+1}} \int_{ \partial B(x_0,R)} \sum_{j=1}^{n} \left(x_0 - y\right)_{j} \left( y-x_0 \right)_{j} f(y)\;d\sigma_{R}(y) \\
&= \lim_{R\rightarrow \infty} \frac{-1}{\omega_{n-1}R^{n-1}} \int_{ \partial B(x_0,R)} f(y)\;d\sigma_{R}(y). \\
\end{split}
\]
Note that for $l \in \mathbb{N}$. On the one hand, we have 
\[
\begin{split}
\int_{l \leq \lvert x_0 - y \rvert \leq 2l} \lvert f(y) \rvert \;dy 
&= \int_{l \leq \lvert x_0 - y \rvert \leq 2l} \frac{\lvert f(y) \rvert}{1 + \lvert y \rvert^{n-2}} \left( 1 + \lvert y \rvert^{n-2}\right) \;dy \\
&\leq l^{n-2} \int_{l \leq \lvert x_0 - y \rvert \leq 2l} \frac{\lvert f(y) \rvert}{1 + \lvert y \rvert^{n-2}} \;dy \leq l^{n-2} \int_{\mathbb{R}^n} \frac{\lvert f(y) \rvert}{1 + \lvert y \rvert^{n-2}} \;dy \;\;\leq l^{n-2},
\end{split}
\]
since $\int_{\mathbb{R}^n} \frac{\lvert f(y) \rvert}{1 + \lvert y \rvert^{n-2}} \;dy$ is finite. On the other hand, by the Spherical Fubini's Theorem and the Mean Value Theorem, we have 
\[
\begin{split}
\int_{l \leq \lvert x_0 - y \rvert \leq 2l} \lvert f(y) \rvert \;dy 
&= \int_{l}^{2l} \int_{\partial B(x_0,R)} \lvert f(y) \rvert \;d\sigma_{R}(y)\;dR = l \int_{\partial B(x_0,R_l)} \lvert f(y) \rvert \;d\sigma_{R_l}(y),
\end{split}
\]
some $R_l \in (l,2l)$. This implies that 
$
\int_{\partial B(x_0,R_l)} \lvert f(y) \rvert \;d\sigma_{R_l}(y) \leq  l^{n-3}.
$
The quantity $I$ becomes 
$
I  = \lim_{l\rightarrow \infty} \frac{-1}{\omega_{n-1}R_l^{n-1}}\int_{ \partial B(x_0,R_l)} f(y)\;d\sigma_{R_l}(y).
$
Therefore, 
\[
\begin{split}
\lvert I \rvert  
&= \limsup_{l\rightarrow \infty} \frac{1}{\omega_{n-1}R_l^{n-1}}\int_{ \partial B(x_0,R_l)} \lvert f(y) \rvert \;d\sigma_{R_l}(y) \leq \limsup_{l\rightarrow \infty} \frac{l^{n-3}}{\omega_{n-1}R_l^{n-1}} \leq \limsup_{l\rightarrow \infty} \frac{1}{\omega_{n-1}l^{2}} = 0.
\end{split}
\]
We now have $I =0$ and $II = f(x_0)$. Hence, $\Delta u(x_0) = f(x_0)$. Since $x_0 \in \mathbb{R}^n$ was arbitrary, this proves that
$
\Delta u = f 
$
in $\mathbb{R}^n$.

We finish this paper by providing the application of the result in the Theorem \ref{mainthm}.
The following theorem shows the existence and uniqueness of the Poisson problem can be proved by assuming the function $f$ belongs to a certain Lorentz space.

\begin{Thm}
Suppose $f \in \mathcal{C}^{1}(\mathbb{R}^n) \cap L^{\frac{n}{2},1}(\mathbb{R}^n)$ such that
$\frac{f(y)}{1+|y|^{n-2}} \in L^1(\mathbb{R}^n)$, and  $\frac{\nabla f(y)}{1+|y|^{n-1}} \in L^1(\mathbb{R}^n)$. Then  the integral
\[
c:= \int_{\mathbb{R}^n} E_n(y)f(y)\;dy
\] 
is well-defined, and the function 
\[
u(x):= \int_{\mathbb{R}^n} E_n(x-y)f(y)\;dy - c, \qquad \text{for}\;\;x \in \mathbb{R}^n
\]
is the unique solution $u \in \mathcal{C}^{2}(\mathbb{R}^n) \cap L^{\infty}(\mathbb{R}^n)$   to the Poisson problem 
\be\label{eqn:P}
  \Delta u = f \quad \text{in} \;\; \mathbb{R}^n\quad\text{and}\quad u(0) = 0.
\ee
\end{Thm}
\noindent {\bf Proof:}
We first claim that $E_n \in L^{\frac{n}{n-2},\infty}(\mathbb{R}^n)$. Let $\lambda \in (0,\infty)$. Then 
\[
\begin{split}
\lbrace  x \in \mathbb{R}^n : \left\vert E_n(x)\right\vert > \lambda    \rbrace  
&= \lbrace  x \in \mathbb{R}^n :  \frac{1}{(n-2)\omega_{n-1}} \frac{1}{\vert x\vert^{n-2}} > \lambda    \rbrace \\
&= \lbrace  x \in \mathbb{R}^n : \left( \frac{1}{(n-2)\omega_{n-1}\lambda}\right)^{\frac{1}{n-2}} > \vert x\vert \rbrace \\
&= B \left( 0,  \left( (n-2)\omega_{n-1}\lambda \right)^{\frac{-1}{n-2}} \right).
\end{split}
\]
Note that
\[
\mu \left(  B \left( 0,  \left( (n-2)\omega_{n-1}\lambda \right)^{\frac{-1}{n-2}} \right) \right)  = c_n \lambda^{\frac{-n}{n-2}},
\]
where $c_n$ is a dimensional constant.
Since $\lambda \in (0, \infty)$ was arbitrary , we have 
\[
\left\Vert E_n \right\Vert_{L^{\frac{n}{n-2},\infty}(\mathbb{R}^n)}=\sup_{\lambda \in (0, \infty)}\left[ \lambda^{\frac{n}{n-2}} \mu \left( \lbrace  x \in \mathbb{R}^n : \left\vert E_n(x)\right\vert > \lambda    \rbrace   \right)\right] = c_n < +\infty, 
\]
which proves the claim. Now, since $E_n \in L^{\frac{n}{n-2}, \infty}(\mathbb{R}^n)$ and $f \in L^{\frac{n}{2},1}(\mathbb{R}^n)$, then by the H{\"o}lder inequality for Lorentz spaces, we have
\[
\left\vert \int_{\mathbb{R}^n} E_n(y)f(y)\;dy \right\vert \leq 
\left\Vert E_n \right\Vert_{L^{\frac{n}{n-2},\infty}(\mathbb{R}^n)} \left\Vert f \right\Vert_{L^{\frac{n}{2},1}(\mathbb{R}^n)} < +\infty.
\]
Thus, the integral $c:= \int_{\mathbb{R}^n} E_n(y)f(y)\;dy$ is well-defined. Next, by the result in Theorem \ref{mainthm}, we have $u \in \mathcal{C}^{2}(\mathbb{R}^n)$, and $\Delta u = f $ pointwise in $\mathbb{R}^n$. Moreover, applying the H{\"o}lder inequality for Lorentz spaces, for any $x \in \mathbb{R}^n$, we have 
\[
\begin{split}
\left\vert u(x) \right\vert 
&\leq \left\Vert E_n(x-\cdot) \right\Vert_{L^{\frac{n}{n-2},\infty}(\mathbb{R}^n)} \left\Vert f \right\Vert_{L^{\frac{n}{2},1}(\mathbb{R}^n)} + \left\vert c \right\vert \leq \left\Vert E_n \right\Vert_{L^{\frac{n}{n-2},\infty}(\mathbb{R}^n)} \left\Vert f \right\Vert_{L^{\frac{n}{2},1}(\mathbb{R}^n)} + \left\vert c \right\vert < +\infty.
\end{split}
\] 
Therefore, $u \in L^{\infty}(\mathbb{R}^n)$. Finally, since $E_n$ is radial, so it is even. Then
\[
u(0) = \int_{\mathbb{R}^n} E_n(-y)f(y)\;dy - c = \int_{\mathbb{R}^n} E_n(y)f(y)\;dy - c = c - c  = 0.
\]
This shows that $u$ is a solution to \eqref{eqn:P}. Suppose $u_1$ and $u_2$ solve the Poisson problem \eqref{eqn:P}. Set $v:= u_1 - u_2$. Then $v \in \mathcal{C}^{2}(\mathbb{R}^n) \cap L^{\infty}(\mathbb{R}^n)$. Furthermore, 
$
\Delta v = \Delta u_1 - \Delta u_2 = 0$ in $\mathbb{R}^n.$
Thus, $v$ is a bounded harmonic function in $\mathbb{R}^n$. By the Liouville's Theorem, $v$ must be constant. Meanwhile, at the origin we have  $v(0)=0$.
Hence, $v$ is identically zero in $\mathbb{R}^n$, which proves the uniqueness. We conclude that $u$ is the unique solution to the Poisson problem \eqref{eqn:P}.

\section*{Acknowledgement}
This work is supported by the SIIT Young Researcher Grant, under a contract
number SIIT 2022-YRG-AS01.  A.C.M.  gratefully acknowledges financial support from the Excellent Foreign Student (EFS) scholarship,  Sirindhorn International Institute of Technology (SIIT),  Thammasat University.
\bibliographystyle{abbrv} 
\bibliography{reference(MS1)}

\end{document}